%% file: AirflowSensing_paper.tex
\documentclass{article}
\usepackage{amsmath,amssymb,amsthm,color}
\usepackage{graphicx}
\usepackage{subcaption}
\usepackage{appendix}
\usepackage[backref,,backend=bibtex,firstinits=true,maxbibnames=9,maxcitenames=6]{biblatex}
\renewbibmacro{in:}{}
\usepackage{caption}
\usepackage{url}
\usepackage{epstopdf}
\usepackage[colorlinks,linkcolor=blue,breaklinks=true,citecolor=blue,pdfstartview=FitH]{hyperref}
\usepackage{todonotes}
\usepackage[font=small,labelfont=bf]{caption}
\theoremstyle{plain}
\newtheorem{theorem}{Theorem}[section]
\newtheorem{prop}[theorem]{Proposition}

\newtheorem{definition}[theorem]{Definition}

\theoremstyle{definition}
\newtheorem{example}[theorem]{Example}

\renewbibmacro*{pageref}{\iflistundef{pageref}{}{\printtext[parens]{\printlist[????pageref][-\value{listtotal}]{pageref}}}}

\theoremstyle{remark}
\newtheorem{remark}[theorem]{Remark}
\usepackage{fullpage}

\title{Sparse sensing and DMD based identification of flow regimes and bifurcations in complex flows}

\date {\today}

\author{B. Kramer\thanks{ Department of Aeronautics and Astronautics, Massachusetts Institute of Technology, Cambridge, MA, USA}\and Piyush Grover  \thanks{Mitsubishi Electric Research Labs, Cambridge, MA, USA}\and Petros Boufounos\footnotemark[2]\and Saleh Nabi\footnotemark[2] \and Mouhacine Benosman\footnotemark[2] }
\graphicspath{{./figures/}}

{\scriptsize
\bibliography{bibliography_boris}}

\newcommand{\opt}{\star}
\newcommand{\conj}{\ast}

\begin{document}

\maketitle

\begin{abstract}
\let\thefootnote\relax\footnotetext{This work was first presented at 2015 SIAM Conference on Applications of Dynamical Systems, Snowbird, UT, USA } 
We present a sparse sensing framework based on Dynamic Mode Decomposition (DMD)
to identify flow regimes and bifurcations in large-scale thermo-fluid systems. 
Motivated by real-time sensing and control of thermal-fluid flows in buildings and equipment, 
we apply this method to a Direct Numerical Simulation (DNS) data set of a 2D laterally
heated cavity. The resulting flow solutions can be divided into several regimes, 
ranging from steady to chaotic flow. The DMD modes and eigenvalues capture the
main temporal and spatial scales in the dynamics belonging to different regimes. 
Our proposed classification method is data-driven, robust w.r.t measurement noise, and exploits the dynamics extracted from the DMD method. 
Namely, we construct an augmented DMD basis, with ``built-in" dynamics, given by the DMD eigenvalues. 
This allows us to employ a short time-series of data 
from sensors, to more robustly classify flow regimes, particularly in the presence of measurement noise. 
We also exploit the incoherence exhibited among the data
generated by different regimes, which persists even if the number of
measurements is small compared to the dimension of the DNS
data. 
The data-driven regime identification algorithm can
enable robust low-order modeling of flows for state estimation and control.

\end{abstract}

\input{sectionIntro}

\input{sectionDMD}

\input{sectionBoussinesq}

\input{sectionClassification}

\input{sectionResults}

\input{sectionconclusions}

\printbibliography

\end{document}

%% file: sectionIntro.tex

\section{Introduction}
The problem of flow sensing and control has received significant
attention in the last two decades. Incorporating airflow dynamics into
the design of control and sensing mechanisms for heating, ventilation
and air conditioning systems yields notable benefits. Combining
such mechanisms with systems that exploit the dynamics of natural convection
to circulate air can save energy and improve comfort. However,
estimation and control of such systems, especially in
real time, is not straightforward. In this paper we provide a framework for sparse sensing and
coarse state estimation of the system.

In principle, the governing equations can be accurately simulated
using DNS techniques. For example, a well accepted mathematical
model for the dynamics of buoyancy driven flows, when temperature
differences are small, is provided by the Boussinesq approximation.
Still, such simulations require vast computational resources, rendering them 
unfeasible in time-critical applications, or use in many-query context.

Moving from simulation to control of flows adds
another level of complexity and poses numerous additional
challenges. In particular, computing a control action based on
full-scale discretized partial differential equation (PDE) models of fluid flow, with typical
dimensions of $n \approx 10^{6-9}$, is computationally prohibitive in
real-time. In addition, in the indoor environments considered herein,
the geometry, boundary conditions and external sources may change over
time. The resulting systems are dependent on large
number of parameters, which further increases their complexity. Of course,
such systems often require sensor-based feedback,
which introduces measurement noise and sensing errors. Thus, modeling and control
methods should be robust to measurement noise and parametric variations.

Undoubtedly, it is of great practical significance to develop a
framework for accurate closed-loop sensing and control strategies for
airflow in a built environment, which are robust to measurement noise and
changing operating conditions. In this paper, we make progress towards
this goal using noise-robust methods based on sparse
detection and classification, which lead to data-driven surrogate models for quick online computation.

For parameter-dependent nonlinear systems, low-order models face additional
challenges. Nonlinear systems can show drastically
different behavior depending on parameters. Modeling strategies which
do not explicitly take this into account, and try to develop `global' parameter independent models instead, are bound to fail. While there have been several
attempts at developing methods to address these issues
\cite{noack14subscaleTerm,noack10ROMgalerkin,sirisup2004spectral,2015arXiv151001728B},
we suggest that it is important to first identify the operating regime, to build a corresponding local low-order model, and only then employ a filtering and
control strategy using these local low-order dynamical models. Since the dynamics of thermo-fluid systems tend to settle on various attractors (such as fixed points, periodic orbits, quasiperiodic orbits), we use the term `regime' to mean `neighborhood of attractors'. This property, shared by many physical systems of interest, has been exploited extensively to develop data-driven multi-scale models \cite{sirisup2005equation,graham1996alternative}.

We propose a method that uses a simple hierarchical strategy. First, for
each possible operating regime of the system, we generate an
appropriate low-order model that captures the short-time spatio-temporal dynamics of the
particular regime. Then, during operation, we use the sensor data to first
detect the appropriate operating regime, the model of which can then be used 
for observing and controlling the system. In this paper, our
main contribution towards this strategy is a general framework for
regime detection, and we also show some numerical results towards coarse reconstruction of the system's state.
Building this framework was possible largely due to the data-driven nature of DMD and the fact that it provides dynamics of subspaces, we heavily exploit in devising our method.

\subsection{Reduced Order Models and Sensing}
The goal of replacing expensive computational models with low-dimensional surrogate models 
in the context of optimal design, control and estimation has led to a rich variety of model reduction strategies in the literature. Nevertheless, there is not yet a
``one-size-fits-all" technique, and each method can outperform others in particular 
applications and settings.

One of the most common methods is Proper
Orthogonal Decomposition (POD)
\cite{aubry88pod,holmes98turbulence,luchtenburg09introPOD,volkwein13notes}. 
This approach finds low-dimensional structures by processing snapshot data from simulation or measurements of 
a dynamical system through a singular value decomposition (SVD). 
The modes are selected by explicitly maximizing the energy preserved in the system. 
%
%
POD-based model reduction has also been used in state estimation of distributed dynamical systems, see
\cite{taylor04PODandLSE,akanyeti11KarmanSensing,bright13CSMLflowCylinder,glauser14airfoilSensing}

Since POD is based snapshots measurements, performance of POD-based can be improved by optimizing the sensors' locations. Using optimized sensor placement, the authors in
\cite{samadiani12datacentersSensors,ghosh14PODdatacenters} employ POD
to predict the temperature profiles in data storage
centers. Furthermore, Willcox~\cite{willcox06gapppyPOD} introduces ``gappy
POD'' for efficient flow reconstruction, and proposes a sensor
selection methodology based on a condition number criterion. Sensor
placement strategies for airflow management based on optimization of observability Gramians
and related system theoretic measures are considered
in~\cite{burns12placement,fang14HVACplacement} using well-established
theory; e.g. \cite{omatu78optimalsensors}, and the references
therein. 

Dynamic Mode Decomposition (DMD)~\cite{rowley09spectralDMD, schmid10DMD,chen12variantsDMD,tu2014dynamic} 
emerged as an alternative to POD for nonlinear systems. This data-driven method 
attempts to capture the dynamics of the system in the low-order model. Therefore, DMD approximates spatial modes and corresponding dynamic information, regarding growth and decay of the modes in time. Computing the DMD requires a few additional steps of, often inexpensive, computation compared to POD.  DMD
has strong connections with the Koopman operator, an approximation of
which is computed in the process~\cite{mezic2005spectral,budivsic2012applied,mezic13analysisKoopman,rowleywilliams15datadrivenDMD}. Various extensions to DMD have also been proposed in the context of sparsity promotion \cite{jovanovic14sparseDMD}, control \cite{brunton16DMDc}, compressed sensing \cite{brunton15CSandDMD}, reduced-order modeling \cite{noack14MORwithDMD}, and 
large datasets \cite{rowleyhemati14DMDstreaming}. 



\subsection{Dynamic Regimes and Classification}
Our work is motivated by recent developments in dynamical systems and
sparse sensing. We only aim to identify what is necessary to
develop an accurate low-order description of the dynamics, and hence,
design a controller. This idea of extracting `effective dynamics' in
complex flows is a central theme in many areas of dynamical systems
and control; see for example \cite{haller2015lagrangian,
  froyland2009almost, ross2012detecting}. Our approach explicitly
takes into account qualitative changes in the dynamics of the system,
as captured by different dynamic regimes. In particular, we formulate
a sensing and detection problem to identify such dynamic regimes. The
analytical and computational methods used herein originate in the sparse
recovery and related literature \cite{eldar10blockRecovery,BKR_TIT11}. 

For parametric systems with different dynamical regimes, the accuracy
of a constructed low-order model depends crucially on choosing the correct
subspace. Hence, in real-time operation, when the system may be in one
of several operating regimes, each corresponding to a distinct
subspace, their classification should be done to identify the relevant modes to be included in a low-order model.

There exists significant literature on the importance of dynamic
regimes. Methods for regime identification, however, are not as
common. In particular, bifurcations in complex systems and their
effect on coherent structures has been studied extensively using the
Perron-Frobenius operator
\cite{gaspard2001liouvillian,junge2004uncertainty,stremler2011topological,grover2012topological}. To
identify bifurcation regimes in a one-dimensional PDE,
\cite{brunton2014compressive} proposes a compressive sensing-based POD
formulation that uses an $\ell_1$ norm optimization to select the relevant elements from a
library of POD-derived bases. In this paper, we build on the framework introduced in \cite{brunton2014compressive} for identification and classification of dynamic regimes through sparse sensing, and extend it in multiple directions.

\subsection{Contributions and Outline}
Our main contribution is a new method that incorporates the dynamic information given by the DMD into 
the sensing and classification step. Using multiple time snapshots from the same sensors,
and thus exploiting the dynamics of the system, we are able to
significantly improve classification accuracy by increasing the
separation between the subspaces.

We demonstrate our general approach on a fluid flow example using the Boussinesq model, which we describe in Section \ref{section:bouss}. This challenging example exhibits several bifurcations as the
governing parameter in the system changes, making it a good
test case to illustrate the proposed method. For this model, the DMD is computed from DNS simulation data
of a 2D laterally heated cavity. 

Then, in Section~\ref{section:classification}, we describe our flow-regime
classification algorithm which combines compressed sensing and sparse
representation techniques~\cite{candes06robust,candes06nearoptimal,donoho06compressedSensing}
with dynamics of the low-order model to classify and reconstruct flow
regimes from few spatial measurements. Fundamentally, our regime classification approach solves a subspace
identification problem, where each regime is represented by a
different subspace. Thus, we develop a theoretical worst-case analysis
that provides classification guarantees. Our guarantees are based on metrics 
from the compressed sensing literature that measure the
separation of different subspaces. The measures can be easily computed for the DMD 
library.

The presented theoretical worst-case classification analysis is conservative, which can be
pessimistic. Our numerical results, described in
Section~\ref{section:results}, demonstrate that the performance is
significantly better in practice. We provide numerical demonstration
of the effectiveness of practical sparse boundary sensing, and a
comparison of the results with a fully distributed array of sensors. We will see that the 
method is particularly robust to measurement noise, and greatly benefits from using the dynamics provided 
from the DMD.

\begin{figure}[h]
\centering
\includegraphics[width=10cm]{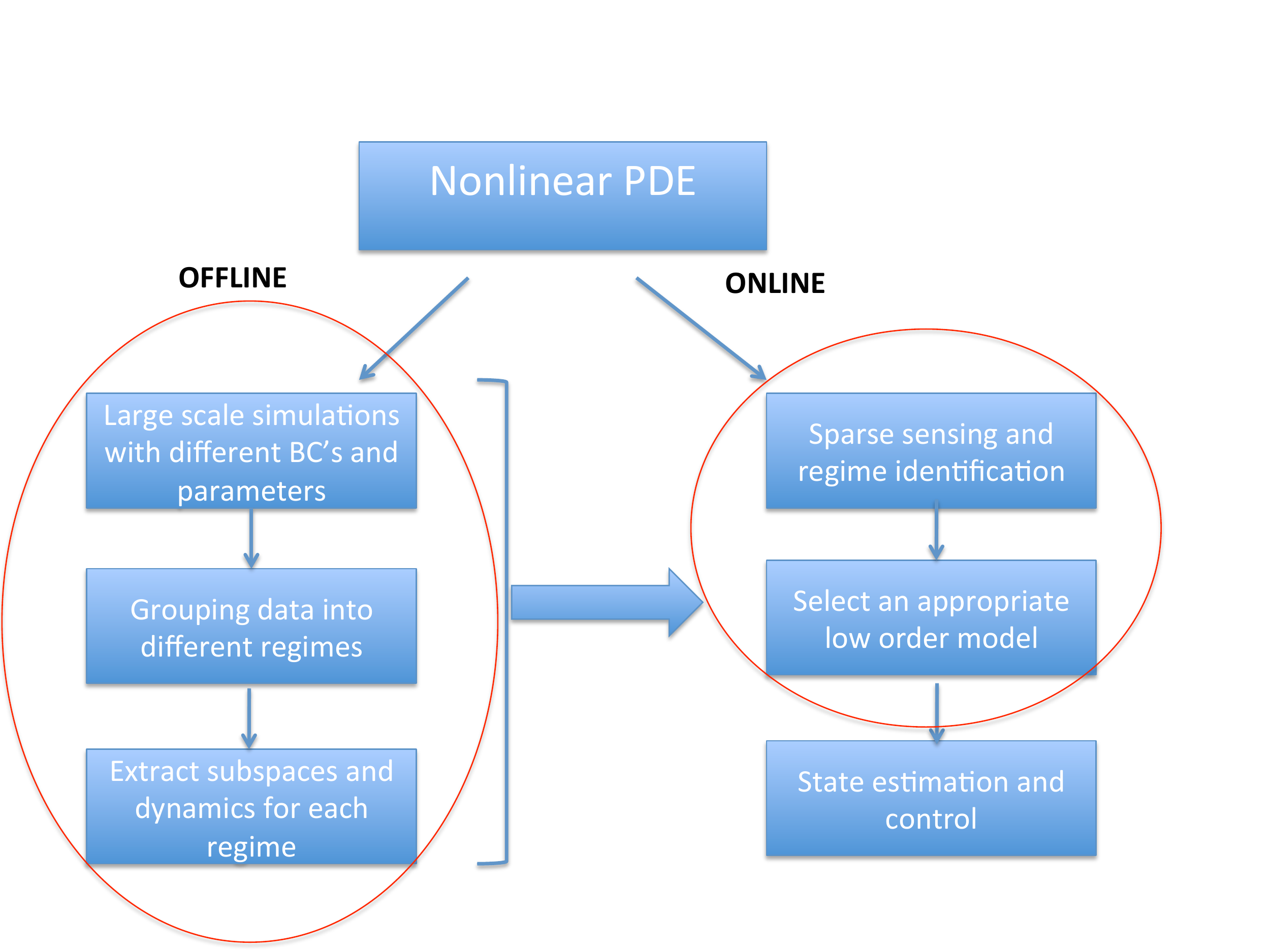}
\caption{\footnotesize{The offline-online approach to model based sensing and
  control of flows, incorporating various operating regimes. In this
  paper, we focus on DMD based offline extraction of subspaces and
  dynamics of various regimes from large scale simulation data, and online sparse sensing based regime
  selection and coarse reconstruction.}}
\label{fig:overview}
\end{figure}

Our overarching goal of developing a closed-loop, low
order flow control system is shown in Figure \ref{fig:overview}. The
work in this paper focuses on the circled components: Large-scale simulation, low-order regime
description, sensing, and data-driven regime identification. Our data-driven framework
provides a solid foundation for developing the whole system.


%% file: sectionDMD.tex
\section{Dynamic Mode Decomposition} \label{section:DMD}
Dynamic mode decomposition is a recently developed data-driven dimensionality
reduction and feature extraction method. It attempts to capture the underlying dynamic evolution of the data, such as that from simulation of a PDE.
Consider a high-dimensional nonlinear system of the form
\begin{equation}\label{eq:highD}
x(t+1)=f(x(t)), \quad x(0) =x_0 \in \mathbb{R}^n
\end{equation}
where our particular interest lies in models (\ref{eq:highD}) arising from spatial discretization of a PDE via, e.g., finite elements, finite volumes, spectral elements. 
Below, we briefly present the DMD formulation following the work of 
in~\cite{schmid10DMD} and~\cite{chen12variantsDMD}. 

Let $s \in \mathbb{N}$ snapshots of the 
high-dimensional dynamical system (\ref{eq:highD}) be
arranged in two $n\times s$ data matrices
\begin{equation*}
X_0 := [x(t_0), \ x(t_1), \ldots, x(t_{s-1})] \qquad \text{and} \qquad X_1 := [x(t_1), \ x(t_2), \ldots, x(t_{s})]. 
\end{equation*}
The dynamic evolution enters into the DMD formulation by assuming that a linear operator $A \in \mathbb{R}^{n\times n}$ maps $X_0$ to $X_1$, namely
\begin{equation}
X_1 = AX_0.     \label{dmdAmap}
\end{equation}
Since the data is finite dimensional, the action of the operator can
be represented as a matrix. The goal of DMD is to approximate the
eigenvalues of $A$, using the data matrices only. 
The first step in DMD is to compute the singular value decomposition of 
$X_0 = W\Sigma V^T$, so that we can approximate the snapshot set via 
\begin{equation}
X_0 \approx X_{0,r} = W_r \Sigma_r V_r^T,  \label{eq:PODstep}
\end{equation}
where $W_r, V_r \in \mathbb{R}^{n\times r}$ contain the first $r$ columns of $W,V$, respectively, and $\Sigma_r$ is the leading $r\times r$ diagonal matrix of $\Sigma$. 
From the Schmidt-Eckart-Young-Mirsky theorem it follows that $||X_0 - X_{0,r}||_2 = \sigma_{r+1}$ (e.g., see \cite[p. 37]{antoulas05approximation}), so if the singular values decay rapidly, the truncation error is small. A scaled version of the left singular vectors $W_r$ are the POD modes of the system. 
In DMD, however, we seek to extract dynamic information about $A$ by considering
\begin{equation}
X_1 \approx A W_r \Sigma_r V_r^T.
\end{equation}
Multiplying by $W_r^T$ from the left and by using the orthogonality of $V_r$, we obtain a reduced-order representation $A_r =  W_r^T A W_r \in \mathbb{R}^{r\times r}$ of the system matrix $A$, namely
\begin{equation}
A_r:= W_r^T X_1 V_r \Sigma_r^{-1},  \label{defSr}
\end{equation}
which is (computationally) much cheaper to analyze than $A$. Next, compute the eigenvalue decomposition
\begin{equation}
A_r Y = Y \Lambda \quad \rightarrow W_r^T A W_r Y = Y \Lambda, \label{eq:DMDSr}
\end{equation}
and by assuming $W_rW_r^T \approx I_r$ one can obtain an approximation for the decomposition 
$$
 AW_r Y \approx W_rY \Lambda,
$$ 
and by defining $\Phi = W_r Y$, it follows that
\begin{equation}
A\Phi \approx  \Phi \Lambda.   \label{dmdAPhi}
\end{equation}
Here,  $\Phi=[\phi_1 \, \phi_2\, \dots\, \phi_r]$ contains the \emph{DMD modes} $\phi_m$  as column vectors. Note, that the matrix $A_r$ is in general non-symmetric, and, therefore DMD modes will be complex. Nevertheless, the $r$-dimensional subspace of DMD modes $\Phi$ in the high-dimensional space $\mathbb{R}^n$ is the same as that spanned by
$W_r$. Hence, the energy content kept in the DMD modes is the same as that for POD modes.

From a computational cost perspective, the dimensionality reduction due to the POD in (\ref{eq:PODstep}) requires a reduced (economy) singular value decomposition of size $n \times s$, where $s \ll n$. The additional computation required in 
the DMD is the eigenvalue decomposition of size $s \times s$ in (\ref{eq:DMDSr}). 

\begin{remark}
	The dynamic mode decomposition provides a non-orthogonal set of modes that attempt to capture the dynamic behavior of the model in a data-driven way. While non-orthogonal bases can burden computations, a breadth of work in DMD has shown that the added benefits of DMD can provide new insights into large-scale models with inherently low-dimensional dynamics \cite{rowley09spectralDMD,schmid10DMD,chen12variantsDMD,tu2014dynamic,mezic2005spectral, budivsic2012applied,mezic13analysisKoopman,rowleywilliams15datadrivenDMD,jovanovic14sparseDMD,brunton16DMDc, brunton15CSandDMD,noack14MORwithDMD,rowleyhemati14DMDstreaming}. We choose to work with DMD for multiple reasons:
	\begin{enumerate}
	\item{
	DMD provides a data-driven \textit{alternative} to other model reduction techniques and, as such, provides great flexibility. Galerkin projection-based model reduction techniques, such as POD, require building a low-dimensional system of ordinary differential equations for future state prediction. This in turn requires having access to, at least, the weak form of the partial differential equation model, and integration routines, which is intrusive. Additionally, problem specific correction terms, such as the shift-modes \cite{noack10ROMgalerkin} or closure models \cite{borggaard14BasisSelection} are often necessary to obtain accurate dynamic information for POD-Galerkin models. }
	\item{
	The information encoded in the DMD modes provide a new viewpoint to the study of low-dimensional, coherent structures in flows, as evidenced in the pioneering work of \cite{rowley09spectralDMD,schmid10DMD}. Every dynamic mode has an associated eigenvalue encoding its dynamic evolution. This yields additional information about spatial structures and their temporal evolution. For our purposes of short-term state prediction, temporal evolution extracted directly via DMD modes and eigenvales is justified via the following arguments. 
Under some conditions on the data
\cite{mezic13analysisKoopman,tu2014dynamic,rowleywilliams15datadrivenDMD}, the DMD modes provide a
linear basis for the evolution of observables, even if the underlying
system is nonlinear. 
In particular, the action of the dynamical system on any finite
set of observables is given by the Koopman operator $\mathcal{K}$,
which is a infinite dimensional linear operator
\cite{mezic2005spectral}. Assuming that the discrete time evolution of
the PDE is given by $x(t+1)=f(x(t))$, the Koopman operator acts on any
scalar observable $g(x)$ as
\begin{equation*}
\mathcal{K}g(x)=g(f(x)).
\end{equation*}
A vector valued observable, such as $h\in \mathbb{R}^n$ can be
expanded in terms of the (scalar) Koopman eigenfunctions of $\mathcal{K}$, $\{\theta_j\}_{j=0}^{\infty}$, and vector valued Koopman modes, $\{\phi_j\}_{j=0}^{\infty}$ as
\begin{equation*}
h(x)=\sum\limits_{j=0}^{\infty} \theta_j(x)\phi_j.
\end{equation*}
Applying the Koopman operator on $h$, we obtain
\begin{equation*}
h(f(x))=\sum\limits_{j=0}^{\infty} \lambda_j\theta_j(x)\phi_j.
\end{equation*}

Given a data set, the unknown scalars $\theta_j(x)$ can be computed,
for example, by projecting on initial conditions, or using some other
way of reducing the overall error in the approximation
\cite{jovanovic14sparseDMD}. In particular, using $\beta_j =
\theta_j(x)$ to simplify notation, the above expansion takes the form
\begin{equation*}
\begin{bmatrix} x(t_1) & x(t_2) & \ldots & x(t_s)\end{bmatrix}\approx \begin{bmatrix}\phi_1 & \phi_2 & \dots & \phi_r\end{bmatrix}\begin{bmatrix} \beta_1 & & & \\ & \beta_2 & & \\ & & \ddots & \\ & & & \beta_r
\end{bmatrix} \begin{bmatrix} 1 & \lambda_1 & \ldots & \lambda_1^{s-1}\\
 1 & \lambda_2 &  \ldots & \lambda_2^{s-1}\\
 \vdots & \vdots & \ddots  & \vdots \\
 1 & \lambda_r & \ldots & \lambda_r^{s-1} \\  \end{bmatrix},
\end{equation*}
where the observable is now a snapshot $x(t_i)$ of the finite-dimensional system \eqref{eq:highD}. In fact, the DMD and the related Koopman mode decomposition have been used to explain the origin and success of various modifications to POD-Galerkin system, such the shift-modes \cite{noack10ROMgalerkin}. 
}
\item{The connection of DMD to many established techniques, such as the Floqu\'{e}t decomposition of linear systems \cite{schmid10DMD}, the Fourier decomposition, and the eigensystem realization algorithm \cite{tu2014dynamic}. Those connections raise further interesting research questions.   }
\end{enumerate}

We note that DMD has several limitations, which are well-documented in the recent literature. While the original DMD algorithm is successful in resolving the dynamics on system attractors (and in some small neighborhoods of these attractors \cite{bagheri2013koopman}), several questions remain on its validity for accurately modeling off-attractor dynamics. However, recent efforts have been directed at accurate extraction of DMD modes (and Koopman modes) for general off-attractor dynamics  \cite{rowleywilliams15datadrivenDMD, mohr2014construction}, and rigorously proving their accuracy in the basin of attraction of these attractors \cite{lan2013linearization}. Hence, as these methods become mature, we expect that our algorithm can be useful for off-attractor dynamics too.

\end{remark}

%% file: sectionBoussinesq.tex
\section{Flow in Two-dimensional Differentially Heated Cavity }\label{section:bouss}

\begin{figure}[h!]
  \centering
    \includegraphics[width=0.5\textwidth]{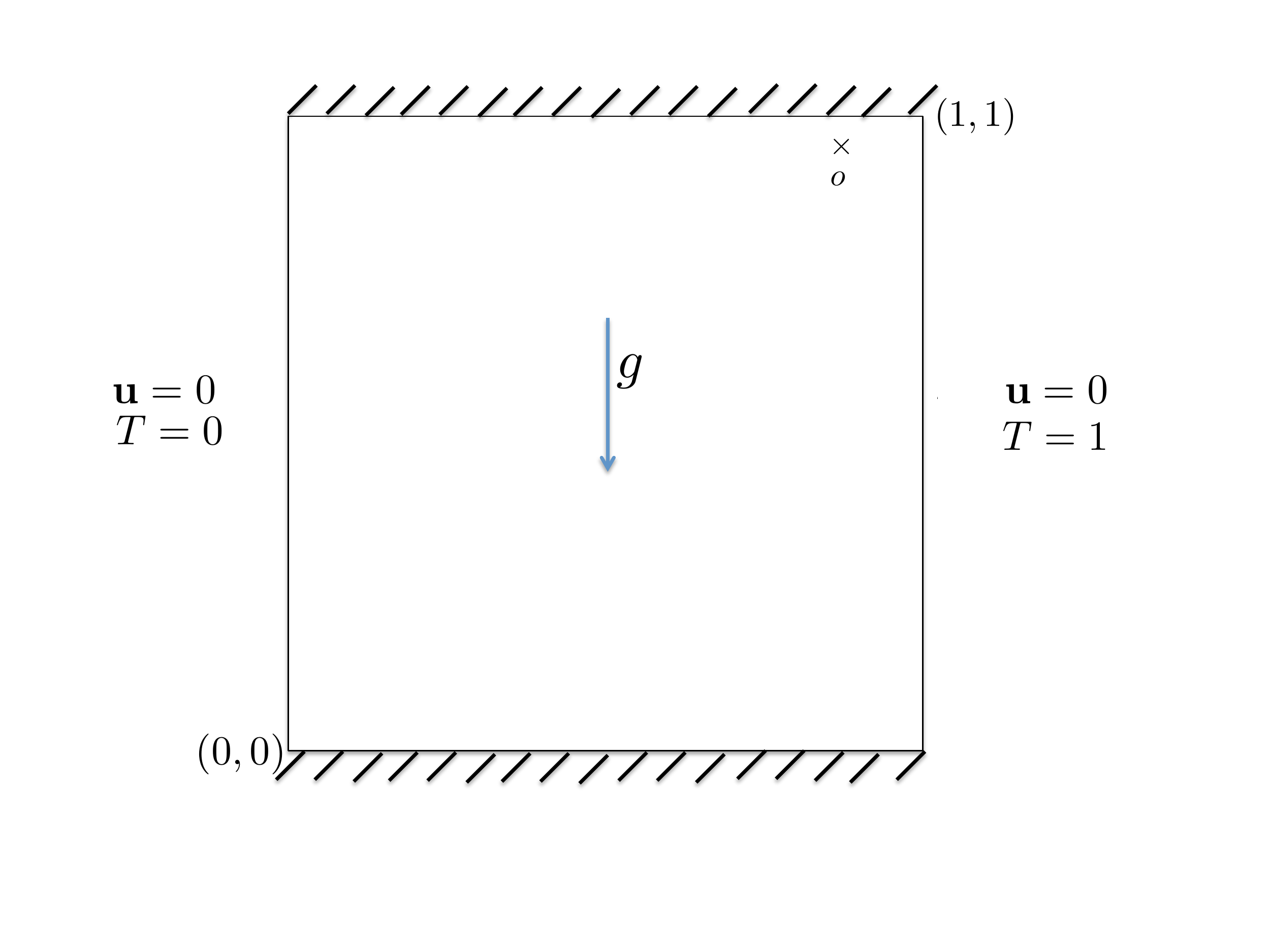}
  \caption{Schematic of the differentially heated cavity problem. The top and bottom walls are insulated and the left and right walls are maintained at a constant distinct
    temperature. The direction of the gravity is downward. We use two monitoring points marked by `o' at $(x,y)=(0.85,0.88)$, and the marker `x' represents location $(x,y)=(0.87,0.93)$.   }
  \label{fig:schem_mesh}
\end{figure}

To demonstrate our approach, we consider the two dimensional upright
differentially heated cavity problem. This is one of the fundamental
flow configurations for heat transfer and fluid mechanics studies, and
has numerous applications including reactor insulation, cooling of
radioactive waste containers, ventilation of rooms and solar energy
collection, among others \cite{Paolucci}.  Figure \ref{fig:schem_mesh}
provides an illustrative schematic of the problem along with the
corresponding boundary conditions.

The domain is a two dimensional enclosure with insulated top and bottom walls, and the left and right walls serve as hot and cold isothermal
sources. We assume the aspect ratio of the cavity is unity;
nonetheless, it should be noted that the variation in aspect ratio
alters the features of the flow non-trivially. The flow is driven due
to buoyancy forces; the temperature difference between the walls
results in a gravitational force exerted on the volume of the fluid
that initiates the flow. The heated fluid rises along the hot wall,
while cooled fluid is falling along the cold wall. When the heated
fluid reaches the top wall, it spreads out to the other side in a form
of the gravity current. In our simulations, the Prandtl number $Pr$ is
0.71, which is a typical value for air. The temperature difference and
other fluid properties are chosen such that the Rayleigh number,
defined as $Ra=\frac{\rho^{2}g \tau \Delta T H^{3}}{\mu^{2}}Pr$,
varies between $10^{1}<Ra<10^{9}$.

\subsection{Governing equations and numerical scheme}
The Boussinesq equations model the viscous, convective fluid motion
associated with buoyancy forces, therefore serving as the model for
natural convection. The Boussinesq equations for an incompressible fluid are given by
\begin{align} 
\nabla \cdot \mathbf{u}&=0  \\
\mathbf{u}_{t} &= \mu \Delta \mathbf{u} - (\mathbf{u} \cdot \nabla) \mathbf{u} - \nabla p + g(\rho-\rho_{0}) \label{eq:governing}\\
T_t &= \kappa\Delta T - \mathbf{u} \cdot \nabla T,
\end{align}
where $T(\cdot,\cdot)$ is the temperature field and
$\mathbf{u}(\cdot,\cdot) = [u_x, u_y]^T$ is the fluid velocity. If the
Boussinesq approximation is applied, the last term in the momentum
equation \eqref{eq:governing} becomes
\begin{equation} \label{eq:approx}
g(\rho-\rho_{0})=\rho_0g\tau(T-T_0),
\end{equation}
where $\tau$ is the coefficient of thermal expansion and $T_{0}$ is
the reference temperature at which $\mu$, $Pr$ and $\rho_0$ are
defined.

To numerically solve the equation and obtain simulation data, we use
the open source spectral-element solver
\texttt{NEK5000}~\cite{Orzsag,patera1984spectral,NEK5000}. Owing to their high
accuracy and general usage, spectral methods are particularly suited
to the study of transition to turbulence in near-to-wall flows, including natural convection within differentially heated
cavities.
\subsection{DNS Results and Discussion}
\begin{figure}[h!]
    \centering
    \includegraphics[height=3in]{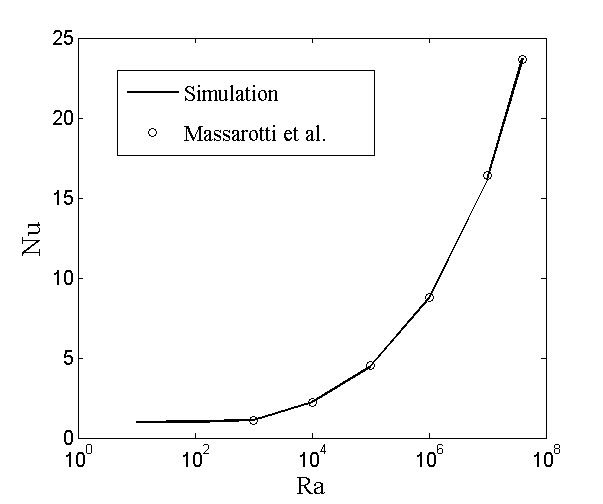}
    \caption{Nusselt numbers from our simulations compared with data from \cite{Massaroti}.}
    \label{figNus}
    \end{figure}

This problem has been extensively studied and accurate solutions are available in the literature for comparison; see, e.g., \cite{Hortmann,Massaroti,quere98diffHeatedCavity}. In Figure \ref{figNus}, we plot the Nusselt numbers computed using our simulation data, along with the corresponding values in \cite{Massaroti}. Our results are in close agreement with the data of \cite{Massaroti}. Since any inaccuracies in resolving near-wall effects will manifest themselves in heat-flux calculations, close agreement in Nusselt numbers with previously reported data in the literature provides evidence of validity of our numerical solver.

A circulatory flow is set up in a vertical layer that is bounded by isothermal
surfaces thermally insulated at the ends, having different
temperatures. The flow ascends against the hot surface and descends at
the cold surface. For smaller Rayleigh number, $Ra<10^{8}$, the flow field reaches a
steady state, that is $\partial \mathbf{u} /\partial t=\partial T/\partial
t=0$. It is well-established that at steady state, the temperature
away from the boundary layers increases linearly over a large part of
the height of the layer~\cite{Paolucci}. Hence, for those values of
Rayleigh numbers, the steady state is stable near the corners of the
cavity as well within the boundary layer close to the walls. 

As the Rayleigh number increases, the flow loses stability. When $Ra\ge Ra_c$,
where $Ra_c$ is critical Rayleigh number, first the flow in the corners
and then the flow in the thermal boundary layer become increasingly
unstable. For large Rayleigh numbers, say, $Ra>10^{9}$ the flow
becomes turbulent. When a statistical steady state is reached, the
space between the vertical boundary layers is filled by a virtually
immobile stably stratified fluid exhibiting low-frequency, low-velocity
oscillations.

\begin{figure}[t!]
\centering
\begin{subfigure}[t]{0.45\textwidth}
\centering
\includegraphics[width=3in, height=2.5in]{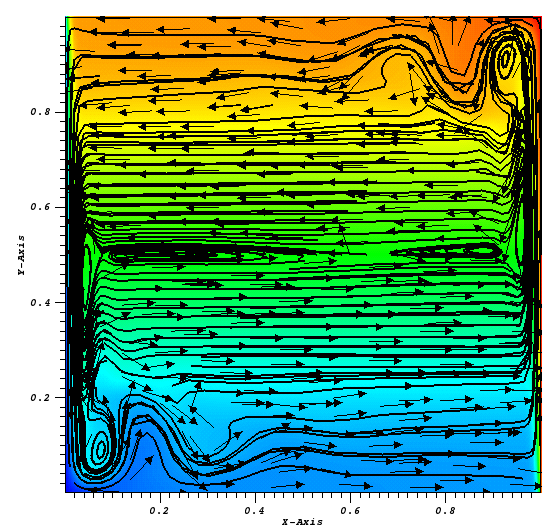}
\subcaption{Streamlines and velocity vectors superposed by temperature psuedocolor (from red to right, the temperature varies from hot to cold).}
\end{subfigure}
\hspace{1cm}
\begin{subfigure}[t]{0.45\textwidth}
\centering
\includegraphics[width=2.5in,height=2.5in]{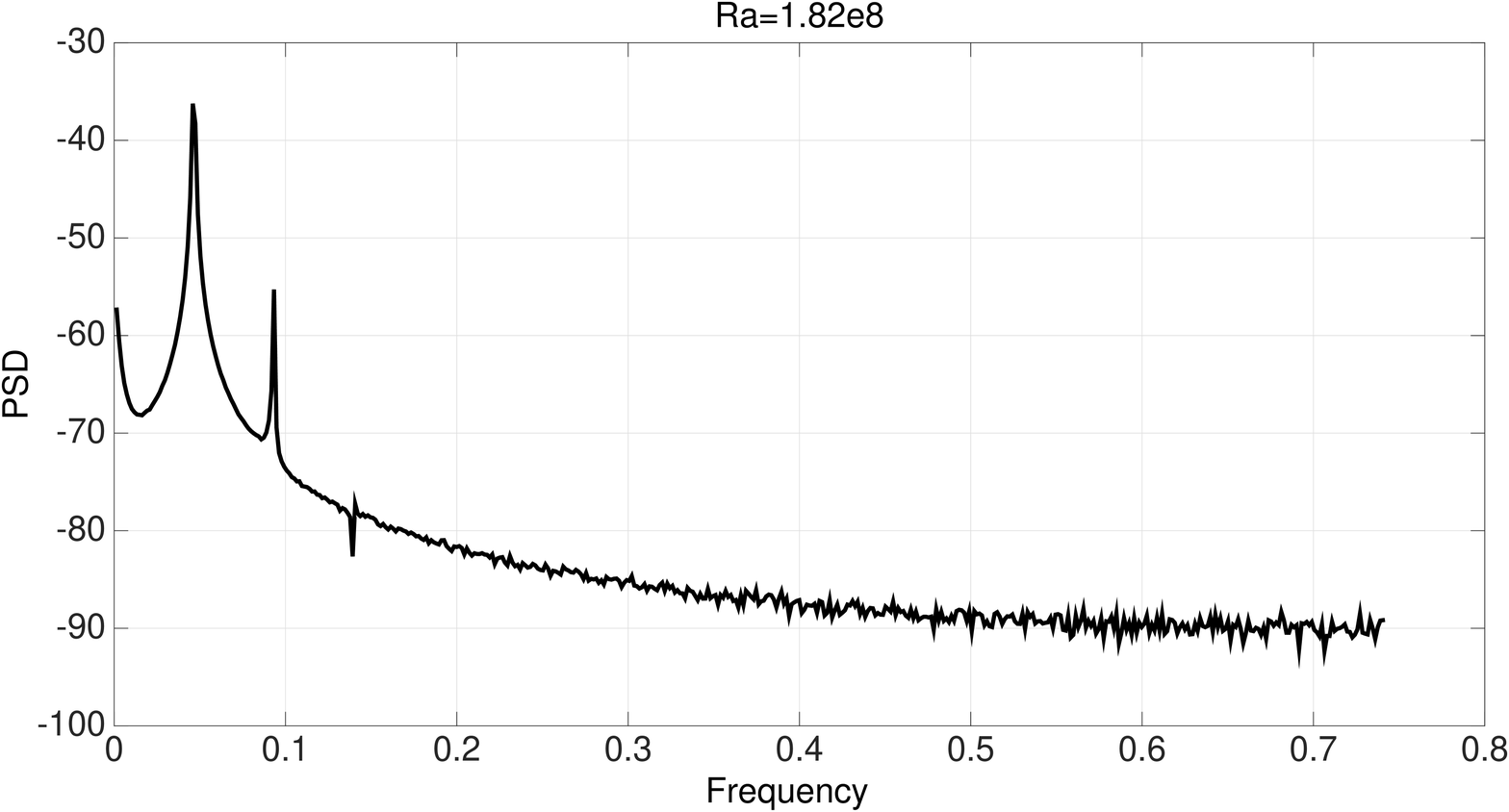}
\subcaption{Power spectrum density of temperature time
	series at location $(x, y) = (0.85, 0.88)$.}
\end{subfigure}
 \caption{Flow solution at $Ra = 1.82\times 10^8$.}
 \label{fig:182e8}
\end{figure}

To better understand the behavior, we examine time series data from local velocity and
temperature measurements at the spatial location $(x,y) =
(0.85, 0.88)$. The temperature at this location reaches a steady-state
asymptotic value after a finite time. By increasing the Rayleigh
number, however, the nature of the flow undergoes a transition.
Consistent with~\cite{quere98diffHeatedCavity}, we observe that the
onset of unsteadiness takes place at a critical $Ra_{c} = 1.82\times
10^8$ and that the first instability mode breaks the usual
centro-symmetry of the solution. Figure~\ref{fig:182e8}(a) illustrates
the streamlines at this Rayleigh number. The onset of unsteadiness in
velocity indicates oscillations in the temperature field as
well. Hence, the time series of the temperature at the monitoring point
shows an asymptotic finite-amplitude periodic state.  The period of
such oscillations is measured from DNS data by using the power
spectral density (PSD) as shown in Figure \ref{fig:182e8}(b).

\begin{figure}[t!]
\centering
\begin{subfigure}[t]{0.51\textwidth}
\centering
    \includegraphics[width=3.5in,height=2.5in]{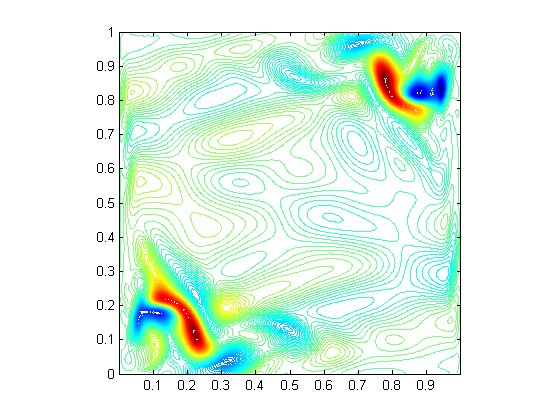}
    \subcaption{Instantaneous fluctuating temperature fields.}
    \end{subfigure}
    \hspace{1cm}
    \begin{subfigure}[t]{0.4\textwidth}
\centering
    \includegraphics[width=3in,height=2.5in]{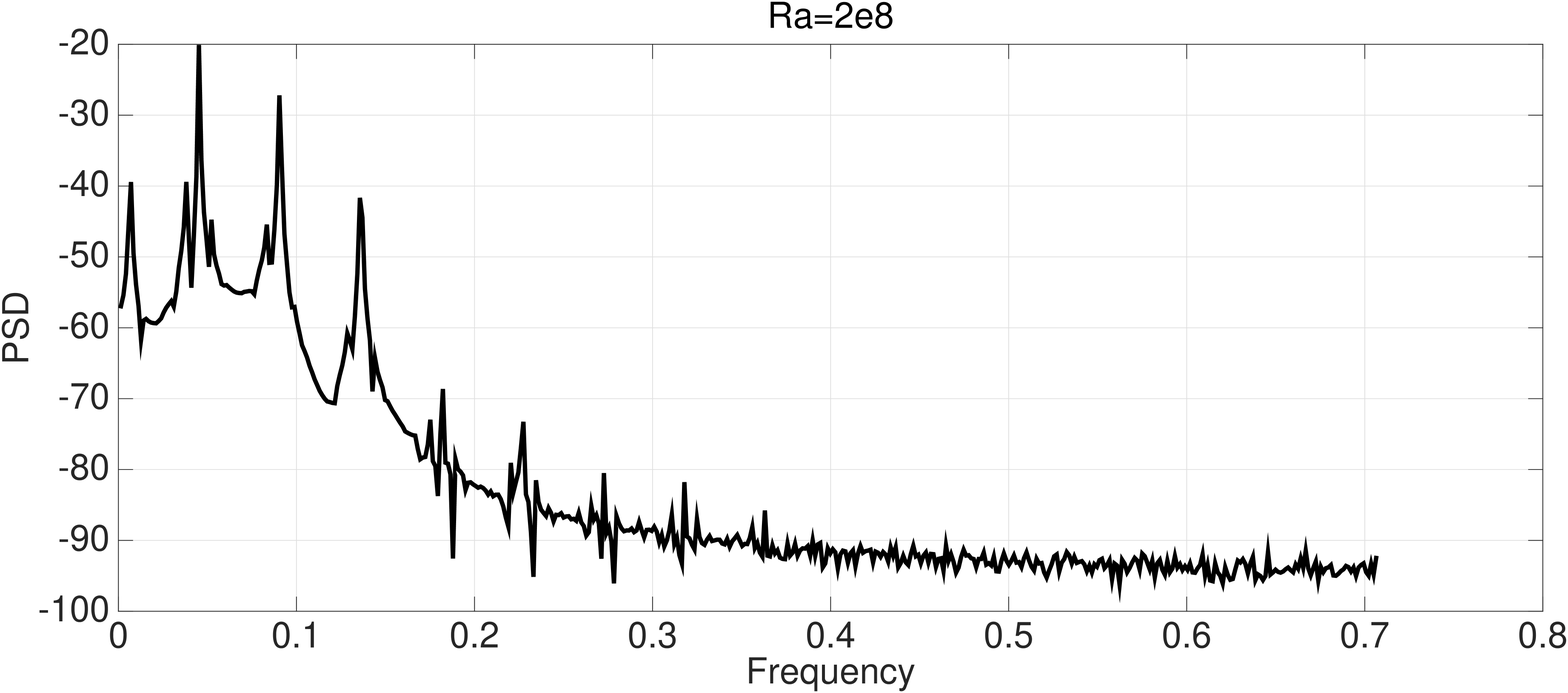}
     \subcaption{Power spectrum density of temperature at
     	location $(x, y) = (0.85, 0.88)$.}
    \end{subfigure}
  \caption{Flow solution at $Ra = 2\times 10^{8}$.}
 \label{fig:2e8}
\end{figure}

This transition to periodic behavior has been studied
in~\cite{quere98diffHeatedCavity}. By numerically computing the
spectrum of the linearized Boussinesq equations operator for Rayleigh
numbers below and above the transition value, it is shown that this
transition is a supercritical Hopf bifurcation.  Physically, this
instability takes place at the base of the detached flow region along
the horizontal walls and is referred to as the `primary instability
mechanism'. The frequency associated with the primary instability mechanism can be
seen in Figures \ref{fig:182e8}(b) and \ref{fig:2e8}(b) to be around
$\omega =0.045Hz$, close to the value reported in~\cite{quere98diffHeatedCavity}. 

With respect to the fluctuating temperature, for $Ra\ge Ra_{c}$, we
observe that away from the corners of the cavity, the contour lines
are inclined at an angle of approximately 20 degrees with respect to
the horizontal, and they propagate in time orthogonally to their
direction. These lines, shown in Figure~\ref{fig:2e8}(a) for $Ra =
2\times 10^{8}$, correspond to the wavefronts of internal waves, which
are shed from the region where the instability mechanism takes place.


\begin{figure}
\centering
\begin{subfigure}[t]{0.46\textwidth}
\centering
        \includegraphics[width=3.5in,height=2.5in]{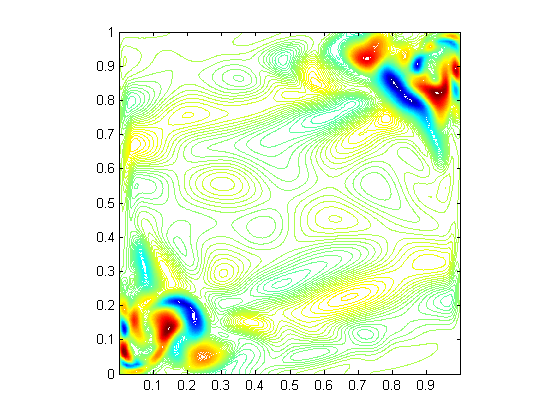}
\subcaption{Instantaneous fluctuating temperature fields.}
\end{subfigure}        
\hspace{1cm}
\begin{subfigure}[t]{0.46\textwidth}
\centering
    \includegraphics[width=3in,height=2.5in]{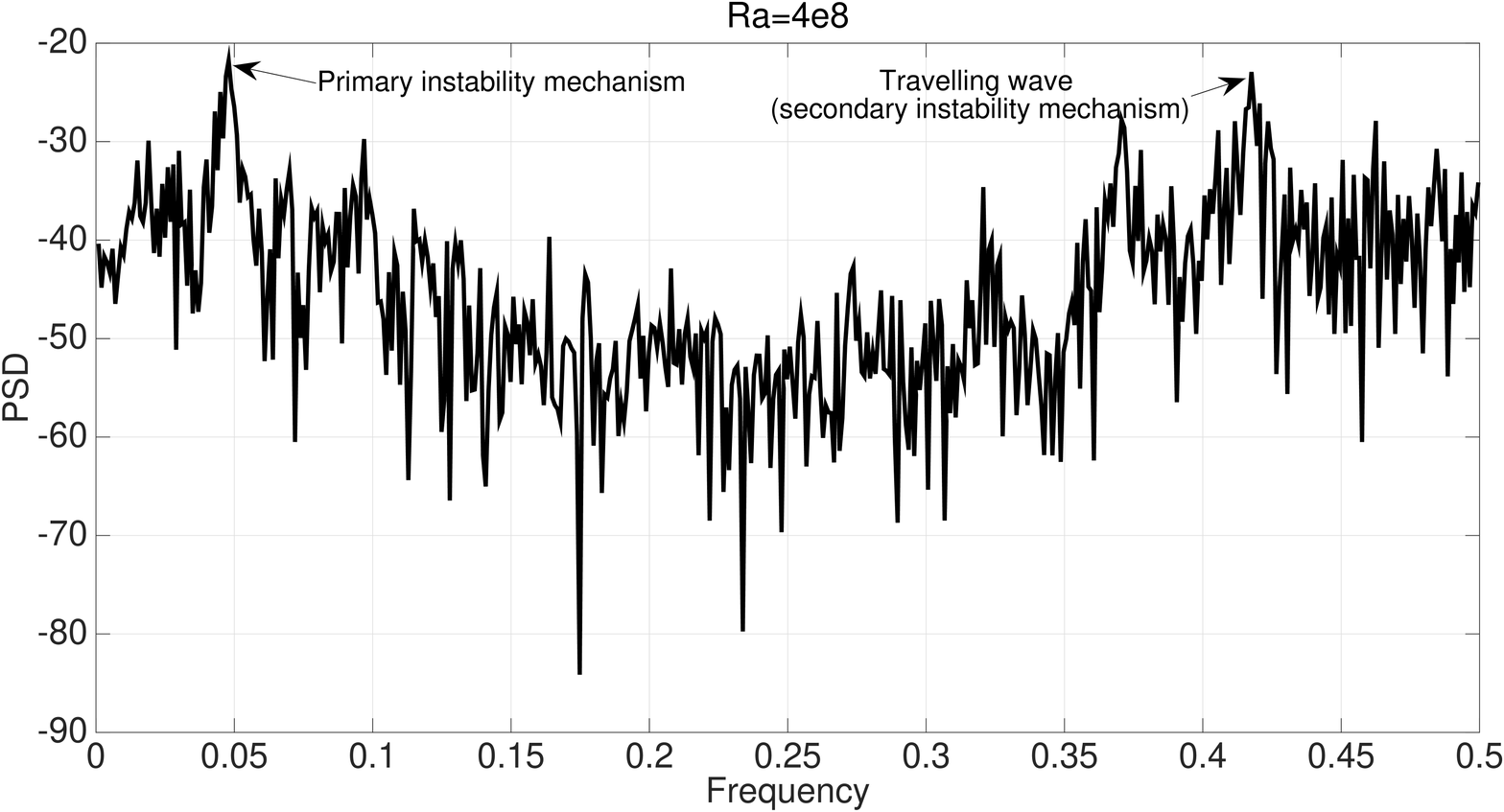}
\subcaption{Power spectrum density of temperature at
	location $(x, y) = (0.87, 0.93)$.}
\end{subfigure}   
  \caption{Flow solution at $Ra = 4\times 10^{8}$.}
 \label{fig:4e8}
\end{figure}

By increasing the Rayleigh number such waves propagate into the domain
as shown in Figure~\ref{fig:4e8}(a) for $Ra = 4\times 10^{8}$. In
addition, we examine the location $(x, y) = (0.87, 0.93)$, closer to
the wall than the one chosen above. At this Rayleigh number, the existence
of a `secondary instability mechanism' has been
reported~\cite{quere98diffHeatedCavity} at a frequency close to
$\omega =0.48Hz$. In Figure~\ref{fig:4e8}(b), we can see that this instability
is also observed in our simulation data. 

The conjecture in \cite{quere98diffHeatedCavity} is that the secondary
instability mechanism is a local phenomenon that originates from the
boundary layer at the top wall. Hence, the monitoring point away
from the wall does not capture the resulting localized oscillations. For even
higher Rayleigh numbers, the solutions depend strongly on the initial condition, and multiplicity of the solutions is also observed. For these higher Rayleigh numbers, the oscillations in
temperature do not show a periodic behavior and the flow becomes
chaotic.

For a variety of parameter settings we compute the DMD from the simulation data. For simplicity we compute the same number
of modes for all different parameters. In particular we find that
$r=20$ modes capture more than $99\%$ of the system energy in all
cases.

\begin{figure}[h!]
\begin{center}
	\begin{subfigure}[t]{\textwidth}
	\hspace{2.1cm}
	\includegraphics[width = 12cm]{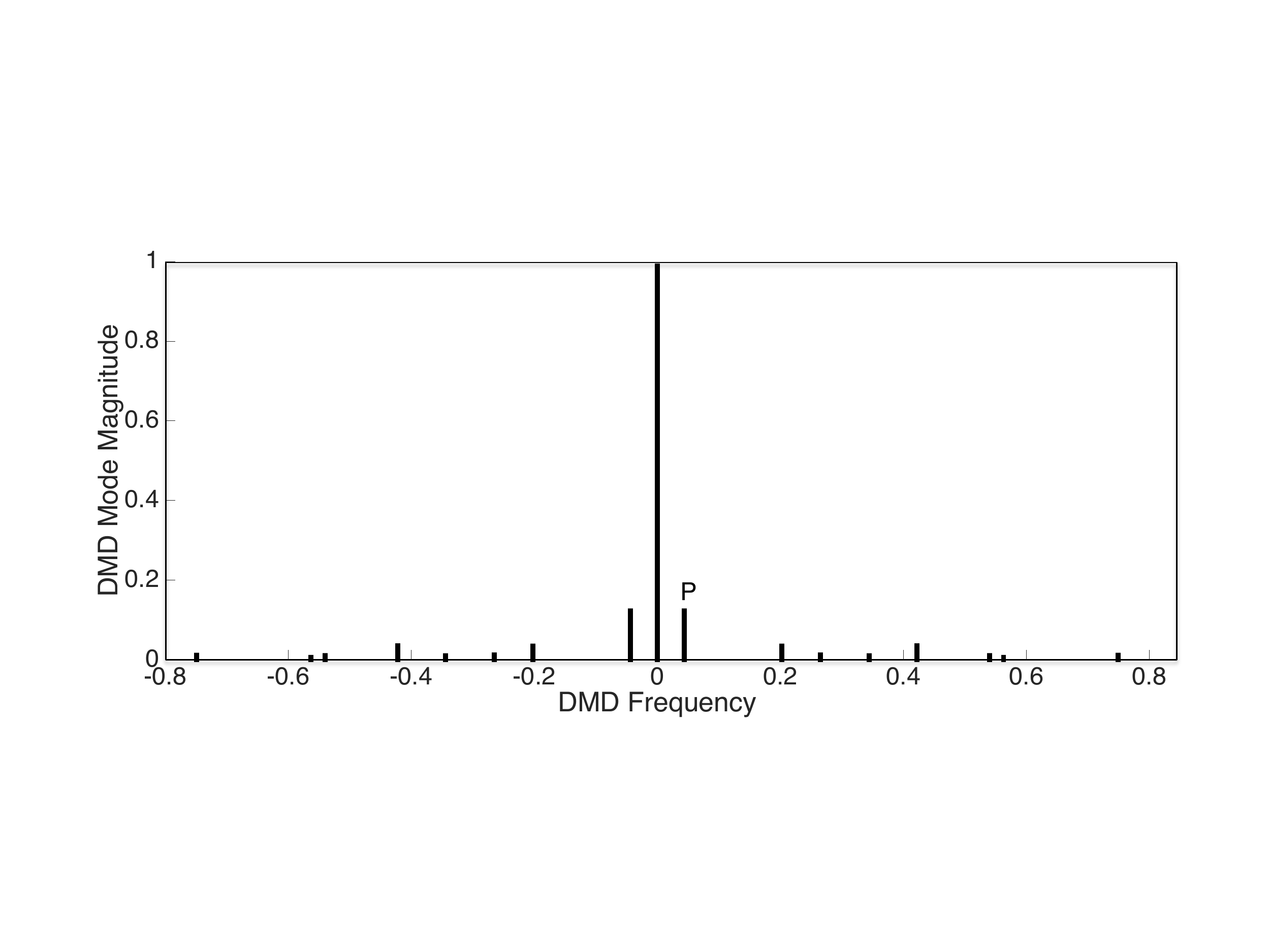}	
	\subcaption{The spectrum of DMD eigenvalues, with corresponding frequencies and DMD mode magnitude. The bar marked 'P' refers to the primary instability frequency discussed in Section \ref{section:bouss}.} 	
	\end{subfigure}
	\hspace{5cm}
	\begin{subfigure}[t]{0.46\textwidth}
	\includegraphics[width = 8cm,height=6cm]{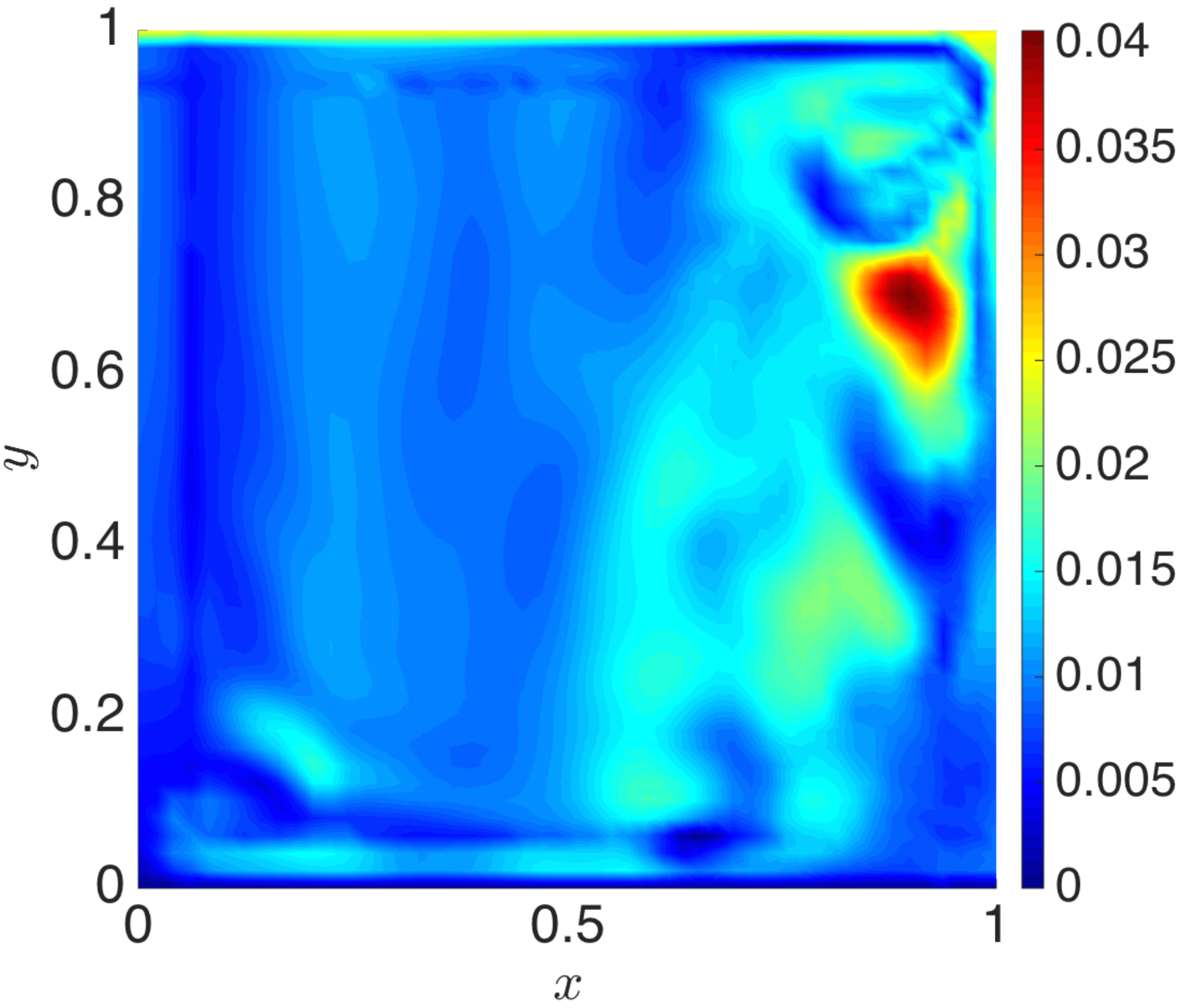}
	\subcaption{DMD mode corresponding to frequency $\omega \approx 0.04Hz$ marked 'P' in plot (a).}	
	\end{subfigure}
	\hspace{1cm}
	\begin{subfigure}[t]{0.46\textwidth}
	\includegraphics[width = 8cm,height=6cm]{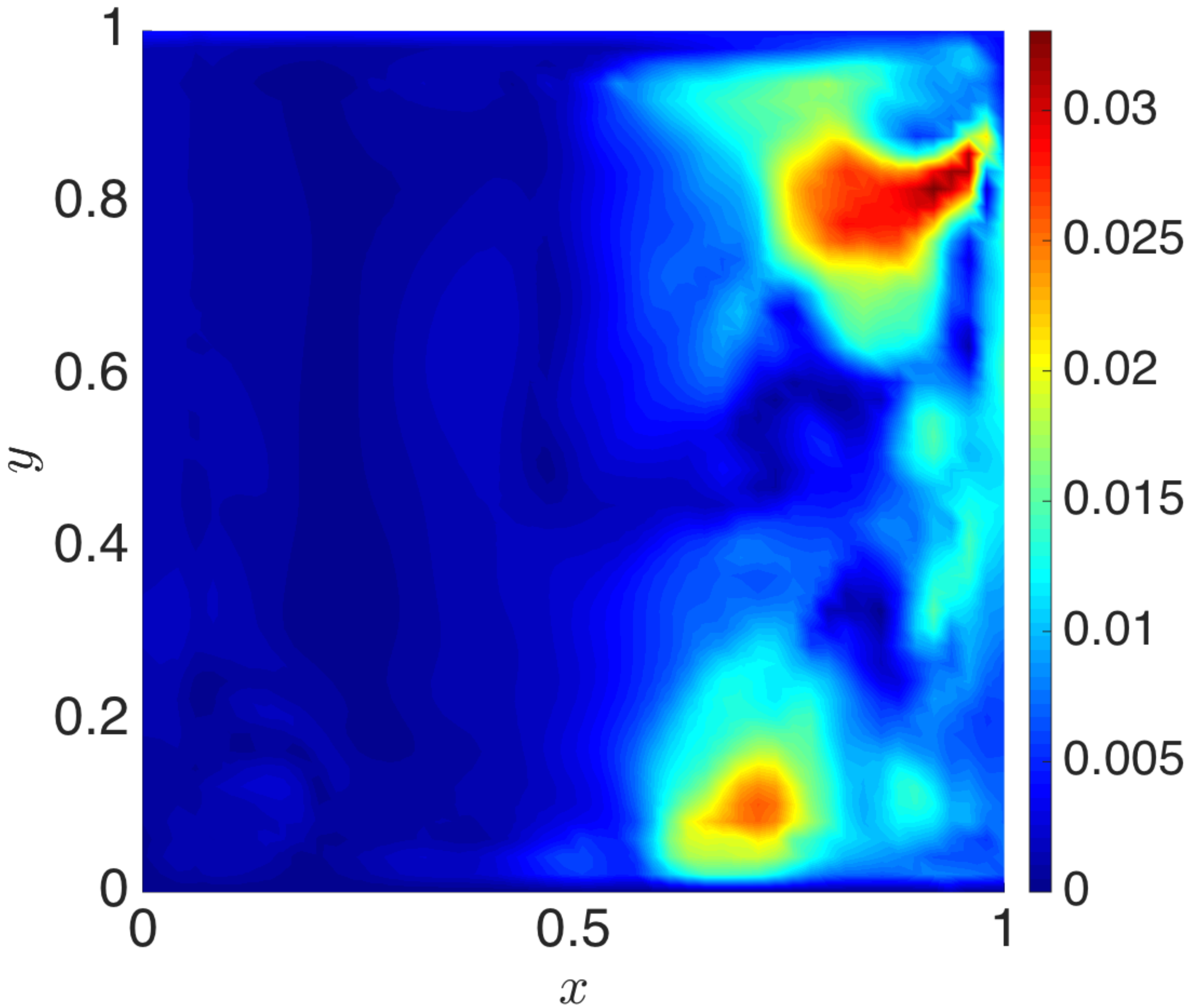}
	\subcaption{DMD mode corresponding to frequency $\omega \approx 0.2Hz$.}		
	\end{subfigure}	
    \caption{DMD modes and spectrum for $Ra=2\times10^8$.}
	\label{fig:DMD_Ra2e8}
\end{center}
\end{figure}

\begin{figure}[h!]
\begin{center}
	\begin{subfigure}[t]{\textwidth}
		\hspace{2cm}
		\includegraphics[width = 12cm]{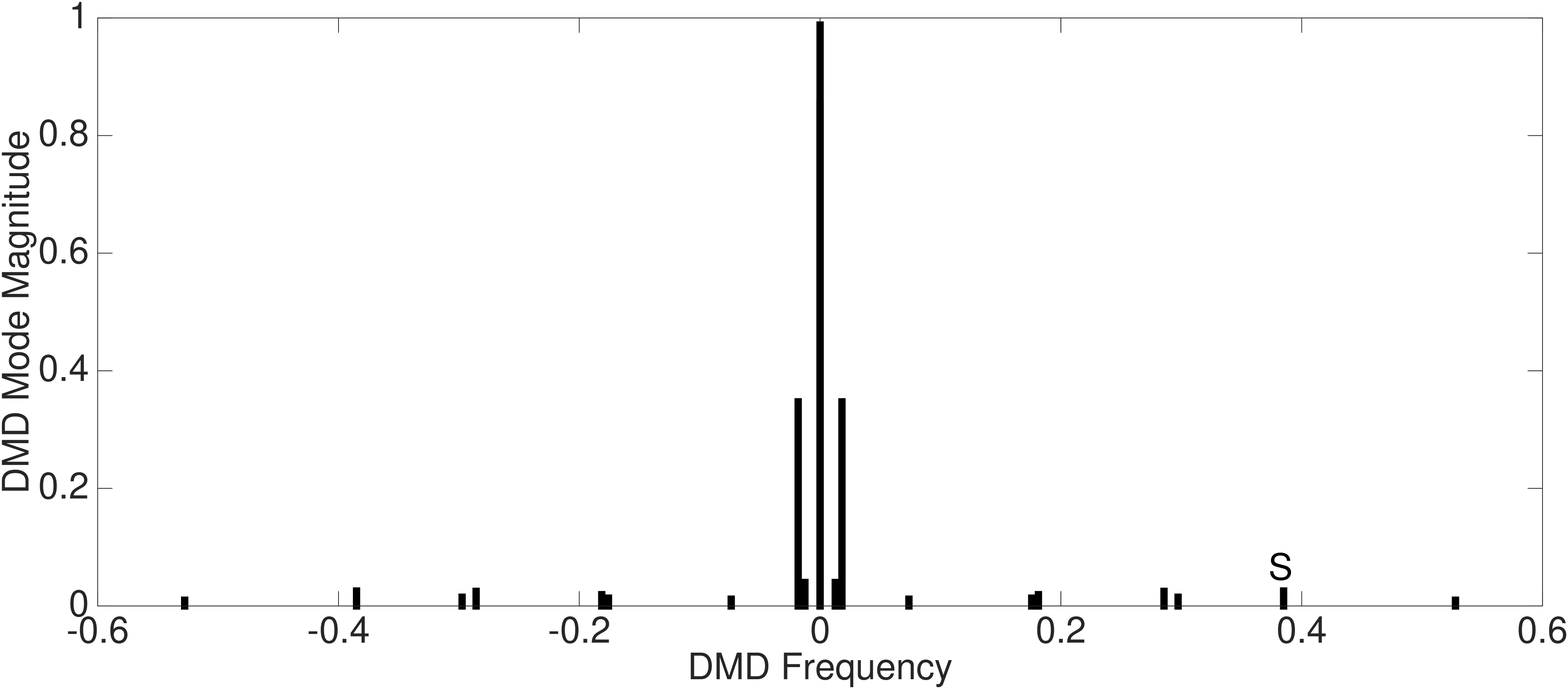}
		\subcaption{The spectrum of DMD eigenvalues, with corresponding frequencies and DMD mode magnitude. The bar marked 'S' refers to the secondary instability frequency discussed in Section \ref{section:bouss}.}
	\end{subfigure}
	
	\begin{subfigure}[t]{0.46\textwidth}
		\includegraphics[width = 8cm,height=6cm]{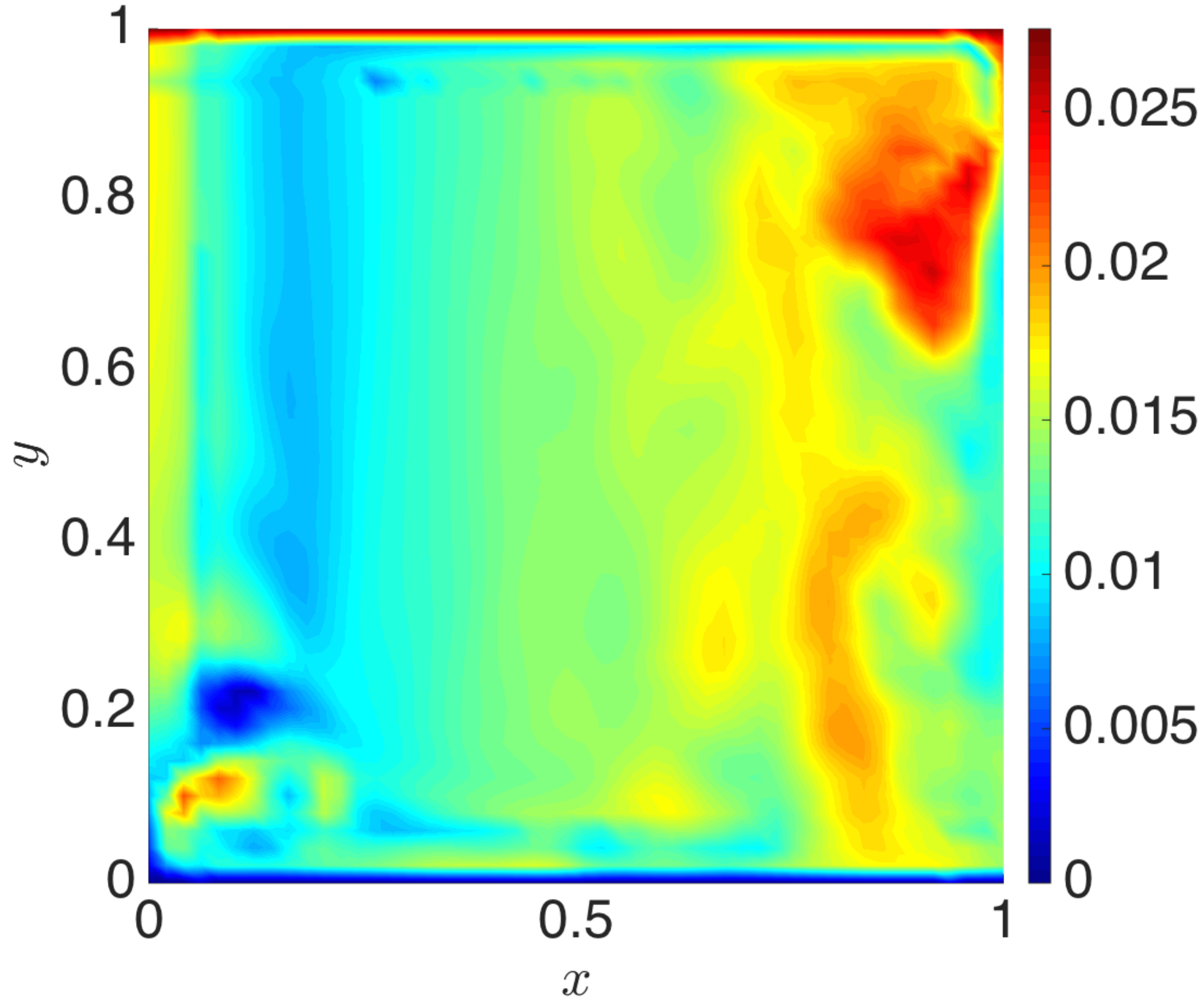}
		\subcaption{DMD mode corresponding to frequency $\omega \approx 0.02Hz$.}
	\end{subfigure}
	\hspace{1cm}
	\begin{subfigure}[t]{0.46\textwidth}
		\includegraphics[width = 8cm,height=6cm]{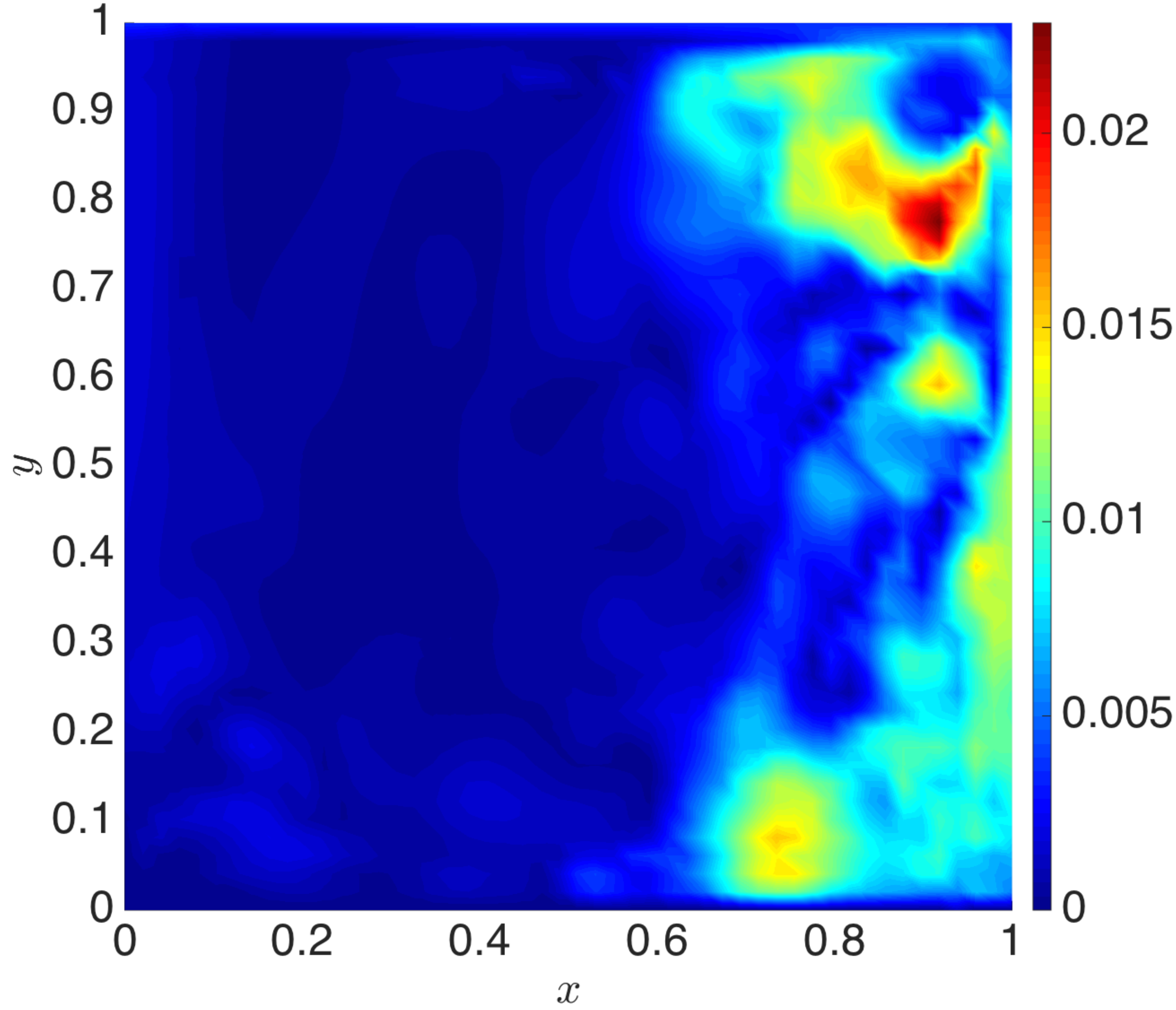}
		\subcaption{DMD mode with frequency marked 'S' $\omega \approx 0.4Hz$.}
	\end{subfigure}			
\end{center}
	\caption{DMD modes and spectrum for $Ra=4\times10^8$. }
  	\label{fig:DMD_Ra4e8}
\end{figure}

We show the dynamic separation that DMD provides for the spatial modes in Figure~\ref{fig:DMD_Ra2e8}(a). We plot the energy content of the DMD modes
at different frequencies obtained from the corresponding DMD
eigenvalues. The dominant DMD mode, other than the `base' flow at zero
frequency, has frequency $\omega \approx 0.045Hz$, close to frequency
associated with the primary instability discussed earlier. The absolute
value of this DMD mode is also plotted over the domain in
Figure~\ref{fig:DMD_Ra2e8}, along with the next dominant DMD mode.
Figure~\ref{fig:DMD_Ra4e8} shows a DMD mode at
frequency $\omega \approx 0.02Hz$ and another DMD mode at frequency close to
the secondary instability frequency $\omega \approx 0.4Hz$. The
mode corresponding to the primary instability was not found among the
dominant DMD modes at this Rayleigh number. This implies that the secondary instability mechanism is not only localized in space, but also has minimal signature in the data in $l_2$ sense.
 For even higher Rayleigh
numbers, DMD modes at several different frequencies have large and
comparable magnitudes, which could be interpreted as a signature of
chaotic behavior.

%

%% file: sectionClassification.tex

\section{Robust Classification by Augmenting DMD Basis}\label{section:classification}

We propose a regime classification approach based on the premise that if a
system operates in a particular regime, then snapshots of the system,
and, therefore of the measurements, lie in a low-dimensional subspace
particular to the regime. 
Thus, our regime detection approach uses offline computation to build a library of subspaces, as introduced in Section \ref{subsec:library}. This general setting allows for any low-dimensional subspace, or reduced-order model technique to be used  for the library generation technique. 
In Section \ref{subsec:class} we introduce the classification problem for a single snapshot in time. We then derive our new method in Section \ref{secAugDMD},  where we augment the DMD basis by using dynamical information. This particular method heavily relies on the dynamic properties extracted from the DMD.
Section \ref{sec:ClassMult} incorporates the augmented DMD basis into the classification of time-sequential measurements of the same dynamic regime, providing a more robust sensing mechanism.
In Sections \ref{subsec:perfAnalysis} and \ref{sec:block-sparse-recov} we then provide a classification analysis, and introduce metrics to assess the quality of the classification.

\subsection{Library generation} \label{subsec:library}
First, we compute a library of dynamic regimes by using subspaces that are generated from a variety of configurations and boundary conditions, each leading to a different regime. Thus, let $\mathcal{Q} =
\{ q_1, q_2, \ldots, q_d \}$ denote the set of $d$ different parameters used to generate the library. For each parameter $q_i \in \mathcal{Q}$, we obtain data from solving the high dimensional model (\ref{eq:highD}). The data are stored in  $X_i := X(q_i) \in \mathbb{R}^{n\times s}$,  where each column of $X(q_i)$ is a snapshot of the solution of (\ref{eq:highD}) associated with parameter $q_i$. Next, we compute $r_i$ basis functions for every regime, $i=1, \ldots, d$, and store them in $\Phi(q_i) := \Phi_i$, a basis for the low-dimensional subspace for the dynamic regime. The set of all such bases $\Phi_i \in \mathbb{R}^{n\times r_i}$ defines a library 
\begin{equation*}
\mathcal{L} := \{ \Phi_1, \Phi_2,  \ldots ,  \Phi_d \}.  \label{defLibrary}
\end{equation*}

\subsection{Classification for single time snapshots} \label{subsec:class}
We use a simple classification algorithm to identify the subspace (and hence the regime) in the library which is most aligned with the current system measurements. 

Our measurements are obtained using a linear measurement matrix $C \in \mathbb{R}^{p \times n}$. The measurement matrix $C$ can represent a large number of possible measurement systems. For example, point measurements of the $j$th component of $x$ are
obtained using rows of $C$ that are zero everywhere except the
$j${th} component, where they take the value of one. Thus, a system
using only distinct single point measurements satisfies
$$C_{i,j} := \{0,1 \}, \qquad \sum_{j=1}^{n} C_{i,j} = 1, \ i=1,\ldots,p, \qquad
\sum_{i=1}^{p} C_{i,j} = 1, \ j = 1, \ldots, n.
$$ 
Alternatively, tomographic measurements, which integrate along one
or more particular directions, can be represented by placing ones on
all the locations along the direction which is integrated. 
Of course, $C$ may also be the identity matrix or another complete
basis for the space. In that case, $p=n$ and the matrix preserves the full
state information from the system. In other words, our regime
identification approach can be used on the full state, if available,
instead of sparse measurements. Since this is not typical in practical
applications, the subsequent exposition assumes that a sensing system
is used, with $p \ll n$.

In practice, the measurements\footnote{For ease of presentation, we start by considering a sensor measurement at a single time instance. We extend this approach by using properties of the DMD in Section \ref{secAugDMD} to use the same spatial sensors, but multiple, time-sequential measurements. As we demonstrate experimentally in Section~\ref{section:results}, the classification estimate then becomes more robust to noise, and leads to high percentages of correct classification. } $y(t) \in \mathbb{R}^p$ are often noisy, so that 
$$
y(t) = Cx(t)+\xi, 
$$ 
where $\xi$ is a noise process, often assumed to be of zero mean and unit variance, i.e., white noise. Assuming that a low-order representation of the $k$th regime represents $x(t)$ in $\Phi_k$, the measurements should lie in the subspace spanned by the $p\times r_k$ matrix
\begin{equation*}
\Theta_k=C\Phi_k. 
\end{equation*}
Correspondingly, a regime library for the observations is given by
\begin{equation}
\mathcal{L}_{obs} = \{ \Theta_1, \ \Theta_2, \ldots, \Theta_d \} . \label{defLibraryTheta}
\end{equation}
Given the measurements, we identify the regime $k^\opt$ as the
subspace in the library $\mathcal{L}_{obs}$ closest to the measurements in an $\ell_2$
sense:
\begin{equation}
k^\opt = \arg \min_{k=1, \ldots,d} \left \{ \min_{\beta_k} \| y(t) -
\Theta_k\beta_k\|_2 \right \},
\end{equation}
where $\beta_k$ denotes the unknown coefficients in the subspace basis
$\Theta_k$. Consequently, $\beta^\opt_k$ is the least square solution to the system
\begin{equation}
\beta^\opt_k = \Theta_k^\dagger y(t),\label{eq:pinv_solution}
\end{equation}
where $(\cdot)^\dagger$ denotes the Moore-Penrose
pseudoinverse. The operator $P_k = \Theta_k \Theta_k^\dagger$ is a
projection operator onto the span of $\Theta_k$ and the classification
algorithm can be expressed as
\begin{equation}
k^\opt = \arg \min_{k=1, \ldots,d} \| y(t) - P_ky(t) \|_2 =
\arg\max_{k=1, \ldots,d}\|P_k y(t) \|_2.
\label{argmax}
\end{equation}
Classification, in other words, projects the measured data $y(t)$ to
the span of each basis set $\Theta_k$ and determines which projection
is closest to the acquired data. \textcolor{black}{Because of the
  orthogonality of the projection error, this is equivalent to maximizing the norm of the projection.}
The underlying classification assumption, which we formalize in the
next section, is that the regimes span sufficiently dissimilar
subspaces and that measurements originating from one regime have
smaller projections onto the subspaces describing the other
regimes. 
If the regime $k$ is identified correctly, then  
$$
x^\opt(t)=\Phi_k\beta^\opt_k=\Phi_k\Theta_k^\dagger y(t), \qquad \Phi_k \in \mathcal{L}, \ \  \Theta_k \in \mathcal{L}_{obs},
$$
approximates the system state $x(t) \in \mathbb{R}^n$ from the measurements $y(t) \in \mathbb{R}^p$. Using $x^\opt(t)$ as an initial state estimate, estimation of the full state vector in future can then be obtained by using the associated local reduced order model (DMD -based or Galerkin approximation-based), and an associated state-space observer, e.g., Luenberger observer for the DMD-based local model.

For the classification given by equation \eqref{argmax} to be successful, the
number of measurements, i.e., the dimensionality of $y(t)$, should be
greater than the dimensionality of the subspace of each regime, i.e.,
$p> \max_k [\mathrm{rank} \ \Theta_k]$. Otherwise the projection error to that
regime is zero. However, the number of measurements needed is significantly less than the $total$ dimension of subspaces of all regimes, $r=\sum_k \mathrm{rank} \ \Theta_k$ (or equivalently, $\sum_k r_k$). In other words, we have the following ordering on the dimension of various quantities, $p\ll r\ll n$, similar to the ordering in \cite{brunton2014compressive}. This ordering is often present in Compressed Sensing literature \cite{candes06CSoverview}. 
This may still imply a large number of sensors if only a single time snapshot is used for classification. This issue is significantly alleviated if we exploit the dynamic information given by the DMD, and use the time augmentation approach we describe in the next section.\\

\begin{remark}
Note that the measurement matrix $C$ might have a nullspace that
eliminates some of the basis elements in $\Phi_i, i=1, \ldots, d$. Subspace-based
identification cannot, therefore, exploit any information along these
basis elements. This can be avoided, for example, using matrices
typical in compressive sensing applications. These include matrices
with entries drawn from i.i.d. Gaussian or Bernoulli distributions,
which are unlikely to have any of the basis elements in their
nullspace---a property known as incoherence (e.g.,
see~\cite{candes06CSoverview}). However, such matrices are
difficult to realize in physical systems. Point and tomographic
sensors are more typical in such systems. Thus, sensor placement, an
issue we briefly discuss later but do not address in this paper, can
significantly affect system performance.
\end{remark}

\subsection{Augmented DMD -- Incorporating dynamics into basis} \label{secAugDMD}
We propose a method to incorporate the dynamic information given by the DMD of the data into the basis generation. 
Consider a state vector $x(t)$
, sampled from the underlying dynamical
system at time $t$. To keep notation minimal, we first consider a single dynamic regime and drop the
subscripts which would indicate regime membership.\footnote{In other words, for this section $\Phi=\Phi_k, r = r_k, x(t) = x_k(t), \beta(t) = \beta_k(t)$.} To begin with, let the state be expressed in the sparse DMD
basis as
$$
x(t) = \Phi \beta(t), 
$$ 
where $\beta(t) \in \mathbb{C}^r$ is again the unknown vector of
coefficients.\footnote{The basis $\Phi$ depends on the data and time
	sampling frequency $\Delta t$. Therefore, in practical sensing, this
	sampling should be kept the same as the one used for generation of
	the basis.} Recall, that by equation
(\ref{dmdAPhi}), the DMD basis approximates the eigenvectors of the
advance operator $A$. Consequently, $A\Phi = \Phi \Lambda$, where
$\Lambda$ denotes the diagonal matrix of the first $r$ eigenvalues of
$A$. From the one-step-advance property of the linear operator
$A$, we have
$$
x(t+1) = Ax(t) = A\Phi \beta(t) = \Phi \Lambda \beta(t), 
$$
and iteratively 
\begin{equation*}
x(t+j) = A^jx(t) = \Phi \Lambda^j \beta(t). 
\end{equation*}
Therefore, subsequent snapshots can be expressed via the \textit{same}
$r$-dimensional vector $\beta(t)$ (and hence the same
regime). Therefore, we have more data to make a confident
classification decision. In particular, the above information can be
written in batch-form as
$$
\begin{bmatrix} x(t) \\ x(t+1) \\ \vdots \\ x(t+j) \end{bmatrix}
= \begin{bmatrix} \Phi \beta(t) \\ \Phi \Lambda \beta(t)
\\ ... \\ \Phi \Lambda^j \beta(t) \end{bmatrix} = \begin{bmatrix}
\phi_1 & \phi_2 & \ldots & \phi_r \\ \lambda_1 \phi_1 & \lambda_2
\phi_2 & \ldots & \lambda_r \phi_r \\ \vdots & \vdots & & \vdots
\\ \lambda_1^j \phi_1 & \lambda_2^j \phi_2 & \ldots & \lambda_r^j
\phi_r \\ \end{bmatrix} \cdot \beta (t) .
$$ 
Next, we define the \textit{augmented DMD basis vector} as
\begin{equation*}
\widehat{\phi}_{i,j} := \begin{bmatrix} \phi_i \\ \lambda_i \phi_i
\\ \vdots \\ \lambda_i^j \phi_i \end{bmatrix} \quad \in
\mathbb{C}^{(j+1)n}, \quad \text{with} \quad \widehat{\phi}_{i,0} := \phi_i,
\end{equation*}
so that the previous equation can be rewritten as
\begin{equation} \label{eq:xaug}
{x}(t:t+j) = \begin{bmatrix} x(t) \\ x(t+1) \\ \vdots
	\\ x(t+j) \end{bmatrix}=
\begin{bmatrix} \widehat{\phi}_{1,k} \ \widehat{\phi}_{2,k}
	\ \ldots \ \widehat{\phi}_{r,k} \end{bmatrix}\cdot
\beta(t) =\widehat{\Phi} \beta (t) .
\end{equation}
When considering the outputs of the dynamical system, $y(t) = Cx(t)$,
the recursion remains unchanged. Thus, using $C \in \mathbb{R}^{p \times r}$ to denote, as above, the sensing matrix, we define
\begin{equation}
\mathcal{C}= \begin{bmatrix} C & 0 &\ldots & 0\\ 0 & C & \ldots &
0\\ \vdots & \vdots & \vdots & \vdots \\0 & 0 & \ldots &
C \end{bmatrix}_{p(j+1)\times r(j+1)}. \label{def:Cblock}
\end{equation} 
and the sensed augmented DMD basis as
\begin{equation}
\widehat{\Theta}_k := \mathcal{C}\widehat{\Phi}_k.
\end{equation}

\noindent We can now define the  library of augmented DMD modes. 

\begin{definition} \label{def_augmented_library}
	Let the DMD modes be $\Phi_i = \Phi(q_i) \in \mathbb{C}^{n\times r_i}$, and let the diagonal matrix of $r_i$ eigenvalues $\Lambda_i = \Lambda_{r_i}(q_i)$ of each dynamic regime $i$, where
	$i=1, \ldots,d$, be given. The {\rm{augmented DMD library}} is defined as
	\begin{equation}
	\widehat{\mathcal{L}}:=\left \{ \begin{matrix} 
	\begin{bmatrix} \Phi_1 \\ \Phi_1 \Lambda_1  \\ \cdots \\ \Phi_1 \Lambda_1^j  \end{bmatrix}, 
	\begin{bmatrix}  \Phi_2 \\ \Phi_2 \Lambda_2 \\ \cdots \\ \Phi_2 \Lambda_2^j \end{bmatrix}, 
	\ldots , 
	\begin{bmatrix} \Phi_d  \\ \Phi_d \Lambda_d \\ \cdots \\ \Phi_d \Lambda_d^j \end{bmatrix}
	\end{matrix}  \right \} 
	=
	\left \{ \widehat{\Phi}_1, \widehat{\Phi}_2, \ldots, \widehat{\Phi}_d  \right \}. \label{defPhiHat}
	\end{equation} 
	Similarly, using $\mathcal{C}$ from \eqref{def:Cblock}, define the observation library as
	\begin{equation*}
	\widehat{\mathcal{L}}_{obs}:= \left \{ \mathcal{C}\widehat{\Phi}_1, \mathcal{C}\widehat{\Phi}_2, \ldots, \mathcal{C}\widehat{\Phi}_d  \right \} = \left \{ \widehat{\Theta}_1, \widehat{\Theta}_2, \ldots, \widehat{\Theta}_d  \right \}.
	\end{equation*}
\end{definition}


\subsection{Classification with time-augmented DMD basis} \label{sec:ClassMult}
To increase robustness of the classification and measuring process, we extend the classification problem from Section \ref{subsec:class} to using multiple time measurements 
$$
y(t:t+j) = [y(t)^T \ y(t+1)^T \ y(t+j)^T]^T.
$$ 
Multiplying \eqref{eq:xaug} by $\mathcal{C}$ from the left (i.e., using only sensed information), the classification problem with the augmented DMD basis can then be recast as finding
\begin{equation*}
k^\opt = \arg \min_{k=1, \ldots,d} \left \{\min_{\beta_k(t)} \| {y}(t:t+j) -
\widehat{\Theta}_k \beta_k(t)\|_2 \right \} .
\end{equation*}
Similarly to the least-squares solution from data at a single time, see ~\eqref{eq:pinv_solution}, the least-squares estimate for the coefficients is computed via the pseudoinverse of
$\widehat{\Theta}_i$, i.e.,
\begin{equation}
\beta^\opt_k(t) = \widehat{\Theta}_k^\dagger y(t:t+j). \label{eq:pinv_solution_augmented}
\end{equation}
In order to cast the classification as a maximization of the
projection energy, as in~\eqref{argmax}, the corresponding projection
operator for $y_k(t:t+j)$ is $\widehat{P}_k =
\widehat{\Theta}_k\widehat{\Theta}_k^\dagger$. Hence, the solution is given by projection,
\begin{align}
k^\opt =\arg\max_{k=1, \ldots,d}\|\widehat{P}_k {y}(t:t+j) \|_2.
\end{align}

Note that $y(t:t+j) \in
\mathbb{R}^{p(j+1)}$ is the available data, and $\beta_i(t) \in \mathbb{C}^{r_i}$. Therefore, we still have
to find $r_i$ coefficients, but this time we can use data of length
$p(j+1)$, where $j$ is the length of the data window to be
specified. In the numerical results reported in Section \ref{section:results}, we will see that $j\leq 10$ is often
sufficient to robustly classify a signal to the correct subspace. 

In
other words, the time augmentation exploits the dynamic information provided by DMD to
increase the amount of data over time, while keeping the number of
spatial sensors the same. This comes at the expense of a small delay waiting to
collect $j$ time snapshots. Moreover, an improvement
is also evident in several library measures that are used in the compressed
sensing community, particularly the alignment and coherence metrics
which we define in the following sections.


\begin{remark}
  Single time-snapshot classification is, of course, possible using
  POD-based reduced models and libraries. However, as evident from our
  development, exploiting the time-evolution dynamics requires the use
  of DMD modes and the DMD-derived eigenvalues. As we see in the
  numerical results below, this significantly increases robustness of
  the classification method, particularly in the presence of sensor
  noise. Moreover, it improves classification accuracy in general.
\end{remark}

\subsection{Classification Performance Analysis} \label{subsec:perfAnalysis}
In this section we develop bounds and metrics that can be used to
guarantee correct classification under the worst-case conditions. As
we observe later in Section~\ref{section:results}, these metrics can
be conservative in practice; however, classification performance is better than
what the bounds suggest. Still, they show the classification
performance trends and provide clear intuition on the role of the
subspaces and their similarity in classification. The metrics we
discuss here and in the next section can be used both with a single
snapshot, using the tools described in Section~\ref{subsec:class}, and
with multiple snapshots in the context of augmented DMD described in
Section~\ref{sec:ClassMult}.

Assume that the system is operating under regime $k$ (we will use
superscripts in this section to denote dependence on the regime) and
that noiseless measurements 
\begin{equation}
y_k(t) = Cx_k(t), \label{ySense}
\end{equation}
are obtained at a single time instance, or over an extended time interval $(t:t+j)$. To
simplify notation, and encompass both cases, for the remainder of this
section we use $y_k$ to denote the measurements, either from a single
or multiple snapshots in time. We also omit $\widehat{\cdot}$ from the
definitions of projections.
The classification algorithm identifies the best matching projection
of the data, $P_{k^\opt}y_k$, and determines the estimated
regime as $k^\opt$. Our goal is to determine metrics and
guarantees under which classification is successful, i.e.,
$k^\opt=k$. To that end we define the subspace $\mathcal{W}_k$
for each $k$, spanned by the basis functions $\Theta_k$.
Given that the measurements originate from regime $k$, most of their
energy lies in $\mathcal{W}_k$. To formalize this statement, we first
decompose the measurements to a direct sum of an in-space
approximation component $\widetilde{y}_k$ and an approximation error
component $\widetilde{y}_k^\perp$, which is orthogonal to the space
\begin{equation}
y_k=\widetilde{y}_k + \widetilde{y}_k^\perp, \qquad \widetilde{y}_k\in
\mathcal{W}_k, \qquad \widetilde{y}_k^\perp\perp
\mathcal{W}_k. \label{xInOut}
\end{equation}
The latter component is due to the approximation performed as part of
the dimensionality reduction. This approximation is accurate if
$\|\widetilde{y}_k(t)^\perp\|_2 \le \epsilon\|\widetilde{y}_k(t)\|_2$ for
some small $\epsilon$ which bounds the approximation error.

Furthermore, we define a metric for subspace alignment, which measures
subspace similarity by determining the vectors in one subspace that
are most similar to their projection in the other subspace:
\begin{align}
\eta_{jk} :=\| P_jP_k \|_2 = \max_{~\forall y}
\frac{||P_j P_k y||_2}{||y||_2}, \quad j,k \in \{1, \ldots,
d\}. \label{defPiPj}
\end{align}  
Using this metric, which is always less than 1, we can show the
following proposition.

\begin{prop} \label{BS_rec_Petros}
Let $d$ subspaces $\mathcal{W}_j, \ j=1, \ldots,d$ be given, and let
the signal $\widetilde{y}_k \in \mathcal{W}_k$ for some $k \in 1, \ldots, d$
according to (\ref{xInOut}), and $t>0$. Moreover, assume that
$\|\widetilde{y}_k^\perp\|_2\le \epsilon\|\widetilde{y}_k\|_2$, with
$\eta_{jk}$ defined in~\eqref{defPiPj}. Then, if
\begin{equation}\label{eq:eta}
 \eta=
\max_{j \neq k} \eta_{jk}<1-\epsilon
\end{equation}
the classification in \eqref{argmax} is successful.
\end{prop}

Before proving the proposition, we provide a brief discussion on the
relevant quantities and a small example application. In particular,
given a set of $d$ subspaces, the $\eta_{jk}$ are quantities easily
computable using simple linear algebra in low
dimensions. Furthermore, the worst-case alignment $\eta$ can be used to formulate the guarantee:
given a set of subspaces, we should expect to always classify signals
correctly if $\eta$ satisfies \eqref{eq:eta}. This in
turn gives a priori guidance on whether a dictionary is suitable for
classification, or if two regimes have similar behavior with
respect to their corresponding subspaces.
In the following example, we provide some intuition on the robustness
of the alignment measure and our bounds.
\begin{example}
  Let us consider two subspaces $\mathcal{W}_j$ and $ \mathcal{W}_k$, and
  assume all signals from regime $k$ contain at least $90\%$ of their
  energy in the subspace $\mathcal{W}_k$, i.e.,
  $\epsilon=.1$. Consequently, if $\eta_{jk} <.9$, we guarantee
  correct classification of the signal $y_k$ to the subspace
  $\mathcal{W}_k$ for all $t>0$.
\end{example}
\begin{proof}
  To demonstrate the proposition, we start with the decompositon
  $y_k=\widetilde{y}_k+\widetilde{y}_k^\perp,~\widetilde{y}_k\in
  \mathcal{W}_k,~\widetilde{y}_k^\perp\perp \mathcal{W}_k,$ as
  above. Since $\widetilde{y}_k^\perp\perp \mathcal{W}_k$, the
  projection onto the correct subspace is equal to
  \begin{equation*}
    P_k y_k = \widetilde{y}_k,
  \end{equation*}
  i.e., has norm equal to
  \begin{equation}
    \|P_k y_k\|_2=\|\widetilde{y}_k\|_2.
    \label{eq:inspace_norm}
  \end{equation}
  The projection onto the other subspaces is equal to
  \begin{equation*}
    P_j y_k = P_j \widetilde{y}_k + P_j\widetilde{y}_k^\perp = P_j P_k
    \widetilde{y}_k + P_j\widetilde{y}_k^\perp,
  \end{equation*}
  which, using the triangle inequality and the trivial bound $\|P_j\widetilde{y}_k^\perp\|_2\le\|\widetilde{y}_k^\perp\|_2$ for any
  projection $P_j$, has norm bounded by
  \begin{equation}
    \|P_j y_k\|_2 \le \|P_j P_k\widetilde{y}_k\|_2 +
    \|P_j\widetilde{y}_k^\perp\|_2 \le \|P_j
    P_k\|_2\|\widetilde{y}_k\|_2 + \|\widetilde{y}_k^\perp\|_2 \le
    \eta\|\widetilde{y}_k\|_2+\|\widetilde{y}_k^\perp\|_2.
    \label{eq:otherspacenorm}
  \end{equation}
  The classification will be accurate if the projection onto the
  subspace corresponding to regime $k$ preserves more energy than all
  other projections, i.e., if
  \begin{equation*}
    \| P_j y_k \|_2 \leq \|P_k y_k ||_2, \quad \mbox{for~all~} j \in
    \{1, \ldots,d \}, j \neq k.
  \end{equation*}
Rewriting equation (\ref{eq:eta}) as
  \begin{align}
    \eta\|\widetilde{y}_k\|_2+\epsilon\|\widetilde{y}_k\|_2<
    \|\widetilde{y}_k\|_2=\|P_k y_k \|_2,
    \label{eq:sufficient_start}
  \end{align}
  where at most a portion $\epsilon$ of the signal lies out of the
  correct subspace $\|\widetilde{y}_k^\perp\|_2\le
  \epsilon\|\widetilde{y}_k\|_2$, i.e.,
  \begin{align}
    \eta\|\widetilde{y}_k\|_2+\|\widetilde{y}_k^\perp\|_2\le
    \eta\|\widetilde{y}_k\|_2+\epsilon\|\widetilde{y}_k\|_2 <
    \|\widetilde{y}_k\|_2=\|P_k y_k \|_2
  \end{align}
  Using~\eqref{eq:otherspacenorm} it follows that for all $j\ne k$
  \begin{align}
    \| P_j y_k \|_2 \leq
    \eta\|\widetilde{y}_k\|_2+\|\widetilde{y}_k^\perp\|_2\le
    \eta\|\widetilde{y}_k\|_2+\epsilon\|\widetilde{y}_k\|_2 <
    \|\widetilde{y}_k\|_2=\|P_k y_k \|_2.
  \end{align}
  Thus, the length of the projection to other subspaces is always
  lower than the length of the projection to the correct subspace and
  the regime is correctly classified. 
\end{proof}

\subsection{Block Sparse Recovery and Classification}  \label{sec:block-sparse-recov}
In addition to the subspace alignment measure $\eta$ introduced in Section \ref{subsec:perfAnalysis}, the
measures of block-coherence between different regimes, and
sub-coherence within each regime, can help to understand the
classification performance and requirements on the
library~\eqref{defLibraryTheta}. These measures are drawn from the
block-sparse recovery and the compressive sensing literature (e.g.,
see \cite{eldar10blockRecovery,BKR_TIT11} and references therein).

Block sparsity models split a vector into blocks of coefficients and
impose that only some of the blocks contain non-zero coefficients.
\begin{definition}
Let $r = \sum_{i=1}^d r_i$ and $\beta_i \in \mathbb{C}^{r_i}$. The vector $\beta = [\beta_1^\conj\ \beta_2^\conj \ \ldots \ \beta_d^\conj]^\conj
\in \mathbb{C}^r$ is called \textit{block $s$-sparse}, if $s$ of its
blocks $\beta(q_i)$ are nonzero.
\end{definition}
Using the above definition, we have that $ y(t) = [\Theta_1, \Theta_2, \ldots, \Theta_d] \ \beta$, with only
$s$ blocks of coefficients in $\beta$ being nonzero, namely the blocks corresponding to the active regimes. We are
interested in conditions on the sensed library $\mathcal{L}_{obs}$,
such that a block-sparse recovery of the vector $x(t)$ from
$p$ measurements is possible. 
\begin{definition} \cite{eldar10blockRecovery}
The \textit{block-coherence} of the library $\mathcal{L}_{obs}$ is defined as 
\begin{equation}
\mu_B := \max_{\substack{i,j=1,\ldots,d \\ i\neq j }} \left [ \frac{1}{r_i} \Vert \Theta_i^* \Theta_j \Vert_2 \right ] , \label{defmuB}
\end{equation}
where it is assumed that $r_1 = r_2 = \ldots = r_d$. 
\end{definition}  
\begin{definition} \cite{eldar10blockRecovery}
The \textit{sub-coherence} of the library is defined as 
\begin{equation*}
\nu := \max_{l\in\{1,\dots,d\}} \max_{\substack{\theta_i, \theta_j \in \ \Theta_l \\ i\neq j} } || \theta_i^* \theta_j ||_2. \label{defSubcoherence}
\end{equation*}
\end{definition}
The sub-coherence gives the worst case measure of non-orthogonality of
various basis elements, computed blockwise. Hence $\nu=0$ if the basis
vectors within each blocks are orthogonal. In particular, when using
DMD, the basis functions are not generally orthogonal, and therefore $\nu \neq 0$. Given the
previous definitions, we can report the following

\begin{theorem} \cite[Thm.3]{eldar10blockRecovery}
\label{thm:block-sparse-recov}
A sufficient condition\footnote{For certain recovery algorithms, such
  as Block-OMP} to recover the block $s$-sparse vector $\beta \in \mathbb{C}^r$ from $y(t) \in
\mathbb{R}^{p}$ measurements via the library $\mathcal{L}_{obs}$ is
$$ s \cdot r < \frac{1}{2} \left ( \frac{1}{\mu_B} + r -
(r-1)\frac{\nu}{\mu_B} \right ),
$$
where it is assumed that all library elements $\Theta_k$ have the same number of
columns, namely $r_1 =r_2=\ldots =r_d$, so that $r=d\cdot r_1$.
\end{theorem}
 Since we are interested in $s=1$ block sparse
solutions (classification of one regime), the above inequality
simplifies to
$$ 
r < \frac{1 +\nu}{\mu_B + \nu} .
$$ 
Note, that the above result provides a sufficient, not a necessary,
condition for accurate recovery of the coefficients $\beta_k$. Thus, the bound in
Theorem~\ref{thm:block-sparse-recov} does not take into account the number
of elements in the block, $r$. Instead, classification corresponds to
recovering only the support, i.e., recovery of the location of the
correct block, a seemingly easier problem. Still, the quantities above
do provide an intuition on the properties of the library of
regimes, including their similarity, in $\mu_B$ and the similarity of
the bases within each regime in $\nu$, affecting the condition number.
When $\nu=0$, i.e., when the basis is
orthonormal, the block-coherence $\mu_B$ is equivalent within a
constant scaling $1/r_1$ to our alignment measure $\eta$.

%% file: sectionResults.tex

\section{Numerical Results}
\label{section:results}
As described in the Section \ref{section:classification}, the data from DNS
simulation is first heuristically divided into 16 regimes. The
spectral element grid for the simulations is finer for higher Rayleigh
numbers, as described in Table \ref{tbl:rayleigh1}. Lagrange polynomials are used as spectral element bases for all simulations.

\begin{table}[h!]
	\centering
	\begin{tabular}{c | c  c  c  c  c c  c  c  c }
		& R1 & R2 & R3 & R4 & R5 & R6 & R7 & R8 & R9 \\ 
		\hline
		$Ra$ & 10 & $10^2$ & $10^3$ & $10^4$& $10^5$ & $10^6$ & $10^7$ & $10^8$ & $1.82\times10^8$ \\ 
		Number of Elements & 4& 4 & 4 & 4& 4 & 4 & 4 & 256 & 256\\
		Polynomial Order & 12& 12 & 12 & 12& 12 & 12 &12 & 18 & 12\\
		State Size $n$ & 1728 & 1728  & 1728   & 1728  & 1728   & 1728   & 1728   & 248,832 & 110,592  \\
		
	\end{tabular}
	
	\vspace{0.5in}
	
	\begin{tabular}{c | c  c  c  c  c  c  c  }
		& R10 & R11 & R12 & R13 & R14 & R15 & R16  \\
		\hline
		$Ra$ & $1.83 \times10^8$  & $1.85\times10^8$ & $2\times10^8$ & $4 \times10^8$ & $6\times10^8$ & $8 \times10^8$  & $10^9$  \\
		Number of Elements &256 & 256& 256 & 256 & 256 & 256 & 256\\
		Polynomial Order &12 & 12& 12 & 12 & 12 & 12 & 18\\
		State Size $n$ &110,592 & 110,592 & 110,592 & 110,592 & 110,592 & 110,592 & 248,832 \\
	\end{tabular}
	\caption{Flow regimes with corresponding Rayleigh numbers and spectral grid specifications.}
	\label{tbl:rayleigh1}
\end{table}


The velocity and temperature data is stacked into the combined state
$x = [u_x \ u_y \ T]^T \in \mathbb{R}^{n}$, and the matrix $X(Ra_i)\in \mathbb{R}^{n\times s}$ contains
the snapshots (in time) as columns. The velocities are
scaled by a factor of 5000, to have both temperature and velocity in
the same order of magnitude. This way, the SVD step of the DMD
 algorithm is not biased towards larger magnitude
entries. Furthemore, the full simulation data is subsequently
interpolated on a $50\times 50$ equidistant grid, so that the data
for all regimes have identical dimension $n=7500$. 

\begin{figure}[h!]
  \centering
  \begin{subfigure}[t]{0.45\textwidth}
  \includegraphics[width=7.45cm]{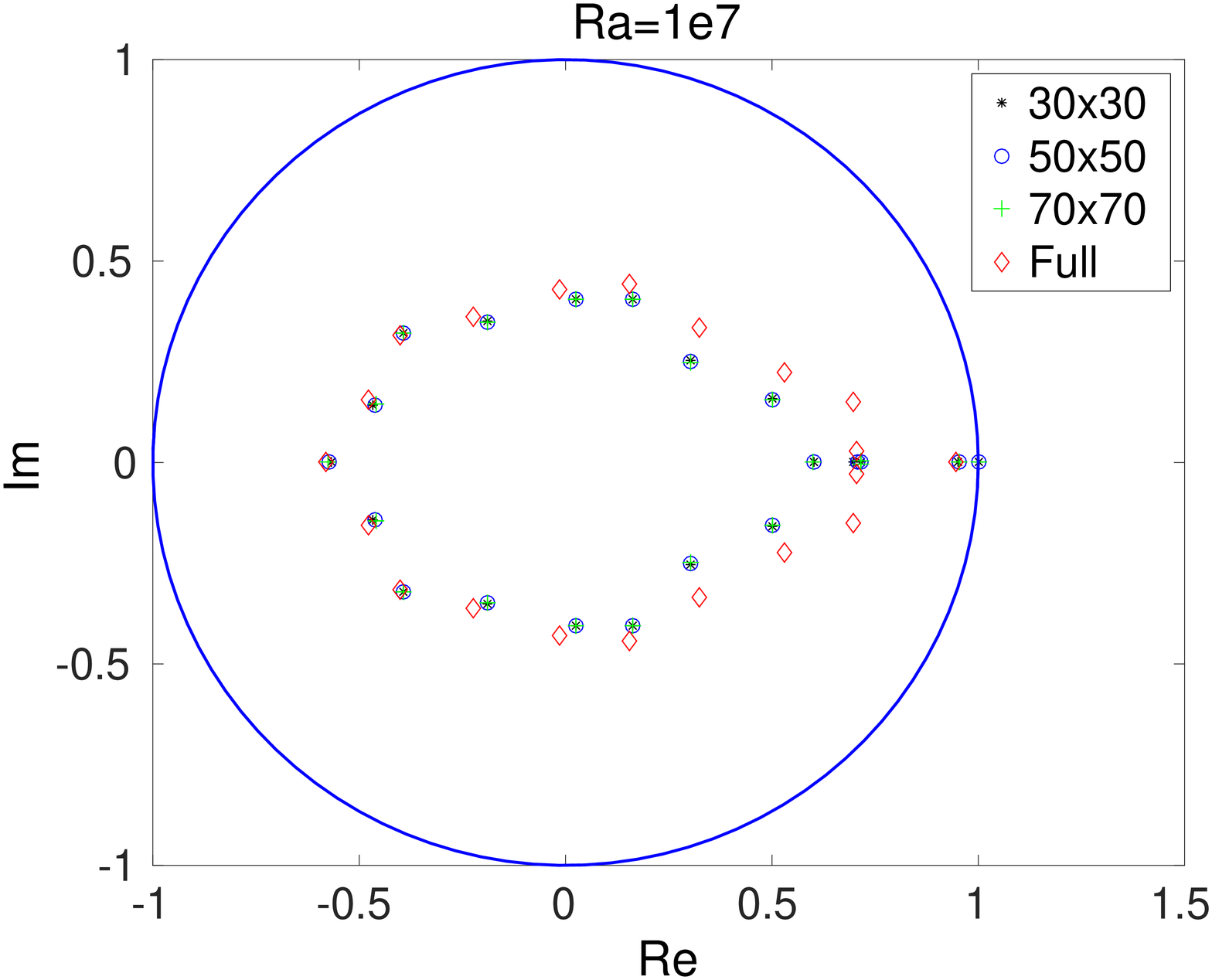}
  \subcaption{ Extrapolated R7 data.}	
  \end{subfigure}
  \begin{subfigure}[t]{0.45\textwidth}
  	\includegraphics[width=7.45cm]{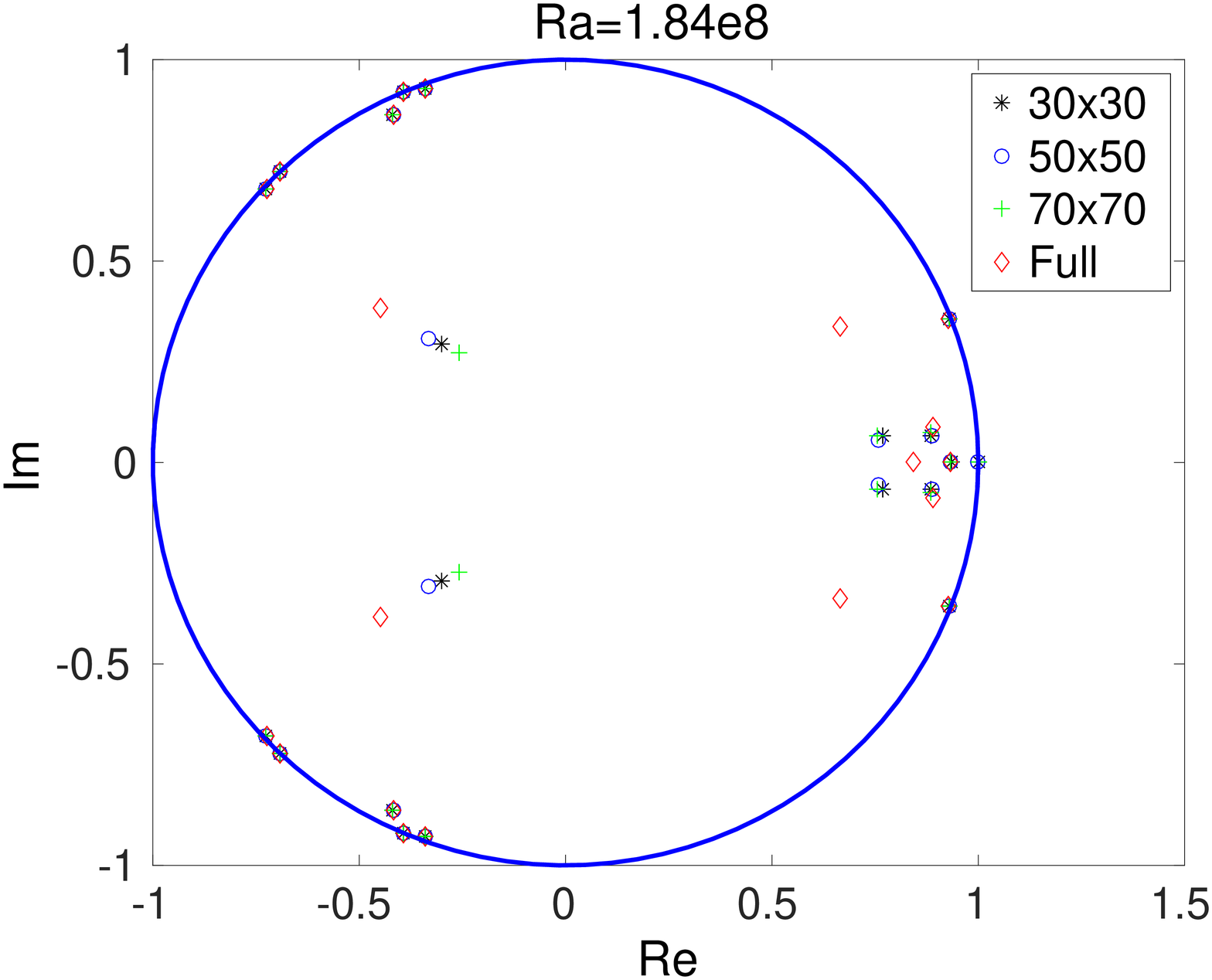}
  \subcaption{Interpolated R11 data.}	
  \end{subfigure}
  \caption{\footnotesize{Comparison of the first twenty DMD eigenvalues computed
    from the full data and interpolated/extrapolated data.}}
  \label{fig:DNS_conv_interp}
\end{figure}

We performed a convergence study with respect to the interpolation
grid size, to ensure that the important information in the flow
solutions is retained. We computed the DMD eigenvalues from the full data 
and compared them with DMD eigenvalues computed from interpolated data.  A good
trade-off between accuracy and size was obtained for the $50 \times 50$ size of the
spectral grid, i.e., $n =7500$. For regimes R1--R7, the solutions are
extrapolated onto this grid, yet this does not change the eigenvalues
considerably. In Figure~\ref{fig:DNS_conv_interp}, a plot of the DMD
spectrum of the first twenty eigenvalues computed from standard DMD is
given for various interpolation sizes. Importantly, the eigenvalues close to the unit circle, which
exhibit mainly oscillatory behavior, converge noticeably quick.

For each of the 16 dynamic regimes given by the parameter in Table \ref{tbl:rayleigh1}, we compute a DMD basis $\Phi_i$ of size $r_i=20$, and subsequently assemble the library of regimes $\mathcal{L}$. 
For the DMD method with augmented basis as introduced in Section \ref{secAugDMD}, we use $10$ additional blocks, corresponding to a sequence of $j=10$ sequential time measurements.
We consider two possible sensing systems and the corresponding
matrices $C$. For both sensing mechanisms, a study with respect to the 
number of sensors is performed, and  the number of sensors used in each test is specified below. In the first sensing system, we sparsely sense from the whole domain by placing point sensors arbitrarily. This serves as a reference, however it is not realistic in
many practical applications to assume that sensors can be placed
arbitrarily. Thus, in the second system, we consider sensing close to the boundary of
the domain, where important boundary effects take place. 

Inspired by compressive sensing principles, both our sensing experiments use randomly selected sensors in the
corresponding sensing area. Thus, we only use a subset of the sensors in the selected area to robustly classify the regimes.

\subsection{Alignment and Coherence Metrics}
We use a variety of metrics to quantify the alignment and coherence of the regime library. The metric $\eta$ from (\ref{defPiPj}) can be interpreted as the worst case similarity measure between two subspaces $\Phi_i$ and $\Phi_j$, since it is based on the spectral norm. We find that for the computed library of regimes, we have $0.95\leq \eta \leq 0.99$. A similar qualitative behavior is seen with $\mu_B$ from \eqref{defmuB}, since it is also based on a spectral norm. Although, as we will see later in this section, $\mu_B$ better captures the decay in coherence among regimes when using augmented DMD.

\begin{figure}[h!]
	\centering
	\begin{subfigure}[t]{0.45\textwidth}
		\includegraphics[width=8cm]{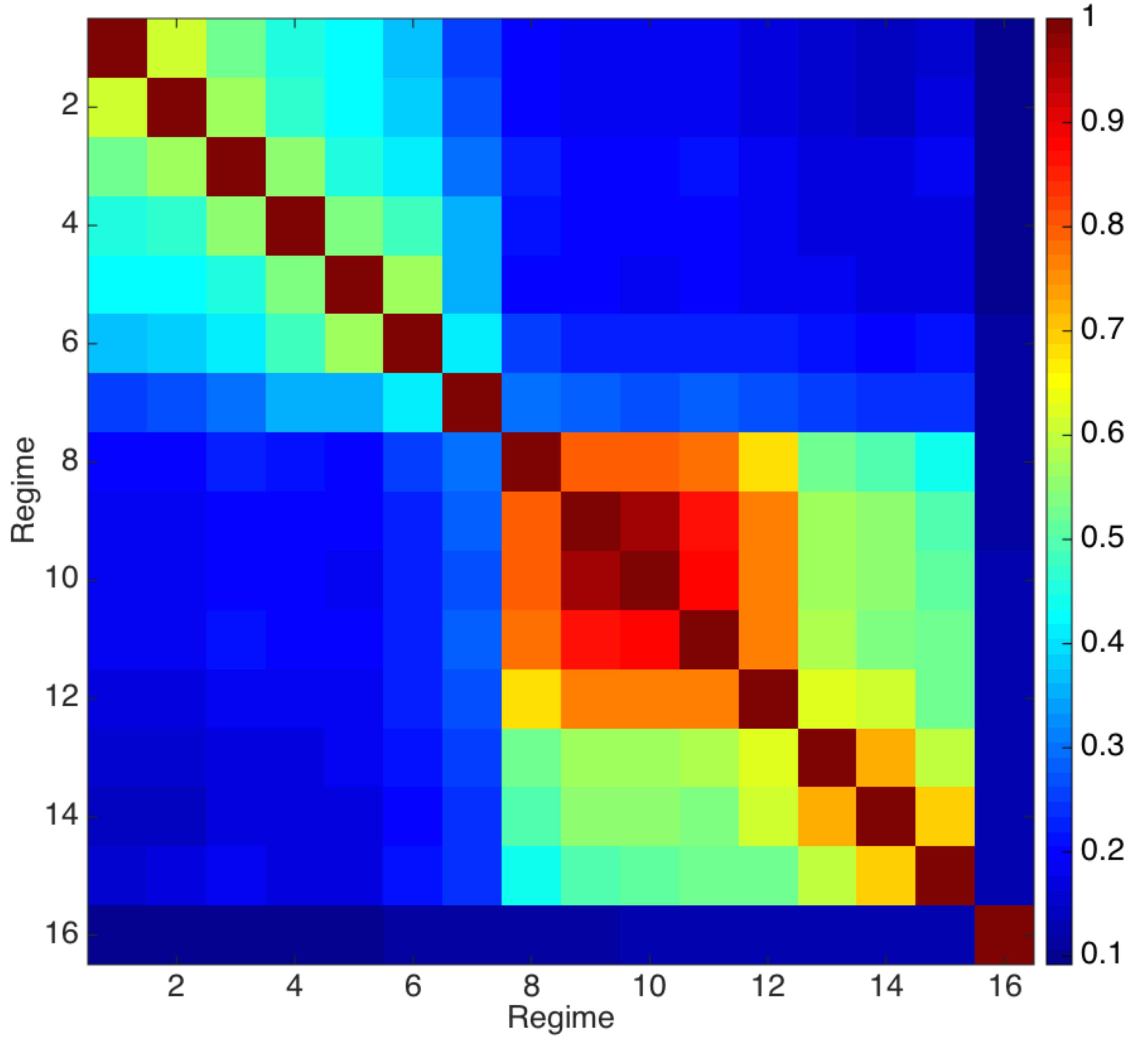}
	\subcaption{Projection onto full $n$-dimensional space, spanned by DMD basis functions, $P_i = \Phi \Phi_i^\dagger$.}	
	\end{subfigure}	
	\hspace{1cm}
	\begin{subfigure}[t]{0.45\textwidth}
		\includegraphics[width=8cm]{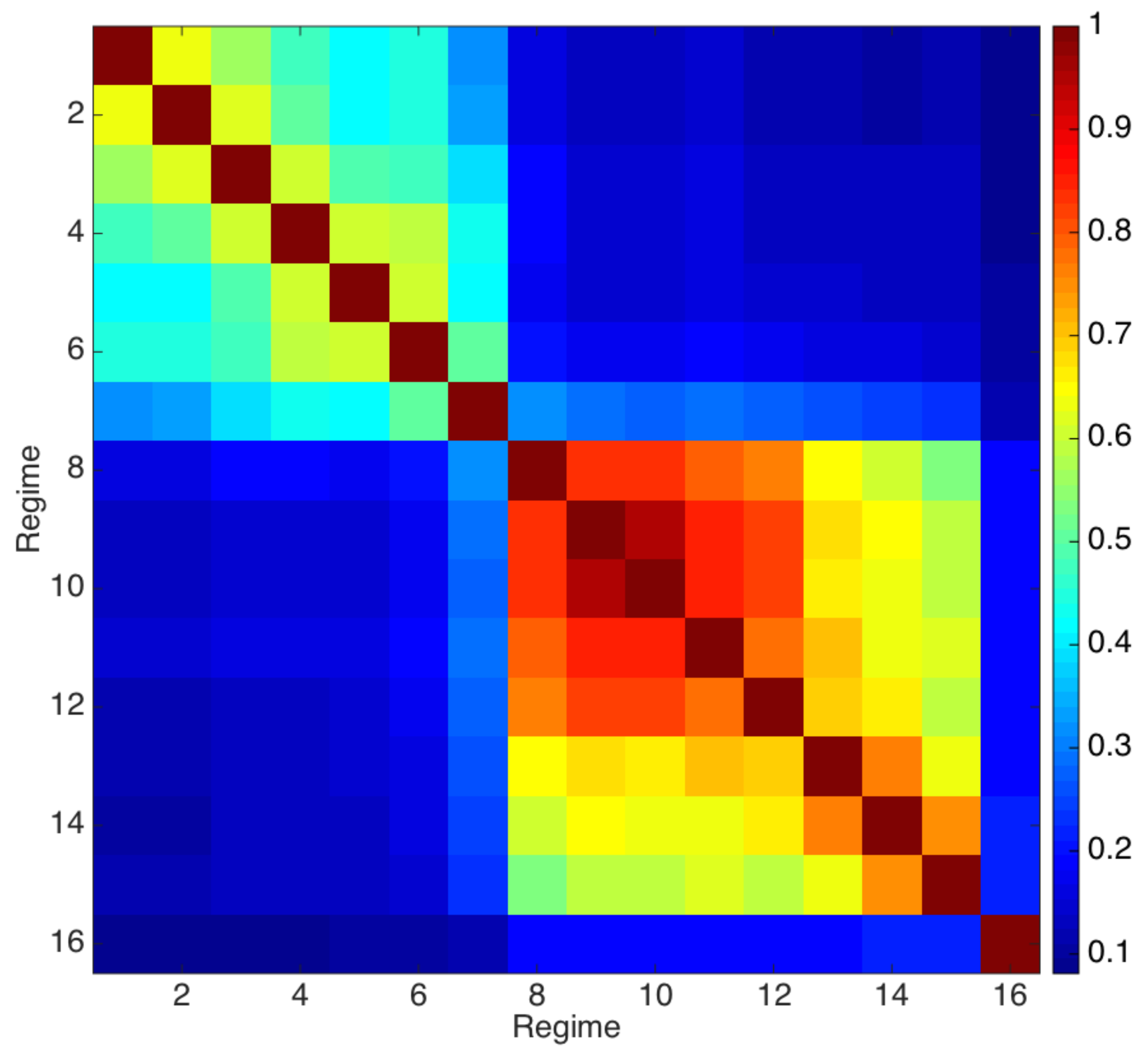}	
		\subcaption{Projection onto the boundary, $P_i=\Theta_i\Theta_i^\dagger$ with $\Theta_i = C\Phi_i$.}	
	\end{subfigure}	
    	\caption{The subspace alignment measure $\gamma_{ij}$ from equation (\ref{def_gamma}).}
	\label{fig:DNS_ssa}
\end{figure}

Hence, to get better insight into the coherence of different regimes, we introduce another measure:
\begin{equation}
	\gamma_{ij}=\frac{||P_i P_j||_F}{||P_i||_F||P_j||_F}, \label{def_gamma}
\end{equation} where $P_i = \Phi
\Phi_i^\dagger$ for the projection of the full state vector, and $P_i=\Theta_i\Theta_i^\dagger$ for the projection subject to measurement $\Theta_i=C\Phi_i$ onto the boundary. The subscript $F$ indicates the Frobenius norm. Figure
\ref{fig:DNS_ssa} shows the measure $\gamma_{ij}$ as defined in
equation (\ref{def_gamma}), both for the full projection, and the projection subject to measurement onto the boundary. 

In contrast to $\eta$, the measure $\gamma_{ij}$ indicates the fraction of information that is retained on \emph{average} by
projecting a random vector onto subspaces of two different regimes
in succession. The diagonal contains ones,
and the off diagonal entries are generally decreasing with the
off-diagonal index, indicating that only neighboring regimes (in terms
of Rayleigh number) share similar features. By definition, the
matrices are symmetric.  Two clusters of regimes appear, the
first one from $Ra=10$ to $Ra=10^7$, and the second cluster from
$Ra=10^8$ to $Ra=8\times10^8$. Note, that the last regime for Rayleigh
number $Ra = 10^9$, resulting in a ``chaotic" flow solution, is
considerably different from the other regimes. Based on this
information, we conjecture that misidentification gets slightly worse
within the two clusters when sensing close to the boundary; we also
expect to see a confusion matrix similar in structure to Figure
\ref{fig:DNS_ssa}. The confusion matrix, used to quantify the success of our algorithm, is defined as follows: the $(i,j)th$ entry contains the percentage of tests in which data from regime $i$ is identified as belonging to regime $j$. Tests are performed independently for data from each regime.

Additionally, we consider the matrix
\begin{equation}
  \kappa_{ij} = \frac{||P_i X_j||_F}{||X_j||_F}, \quad i,j = 1,
  \ldots, d, \label{def_kappa}
\end{equation}
which measures the energy in the projected subspace compared to the
actual data, and allows estimation of $\epsilon$.
 
\begin{figure}[h!]
	\centering
	
	\begin{subfigure}[t]{0.45\textwidth}
	\includegraphics[width=8cm]{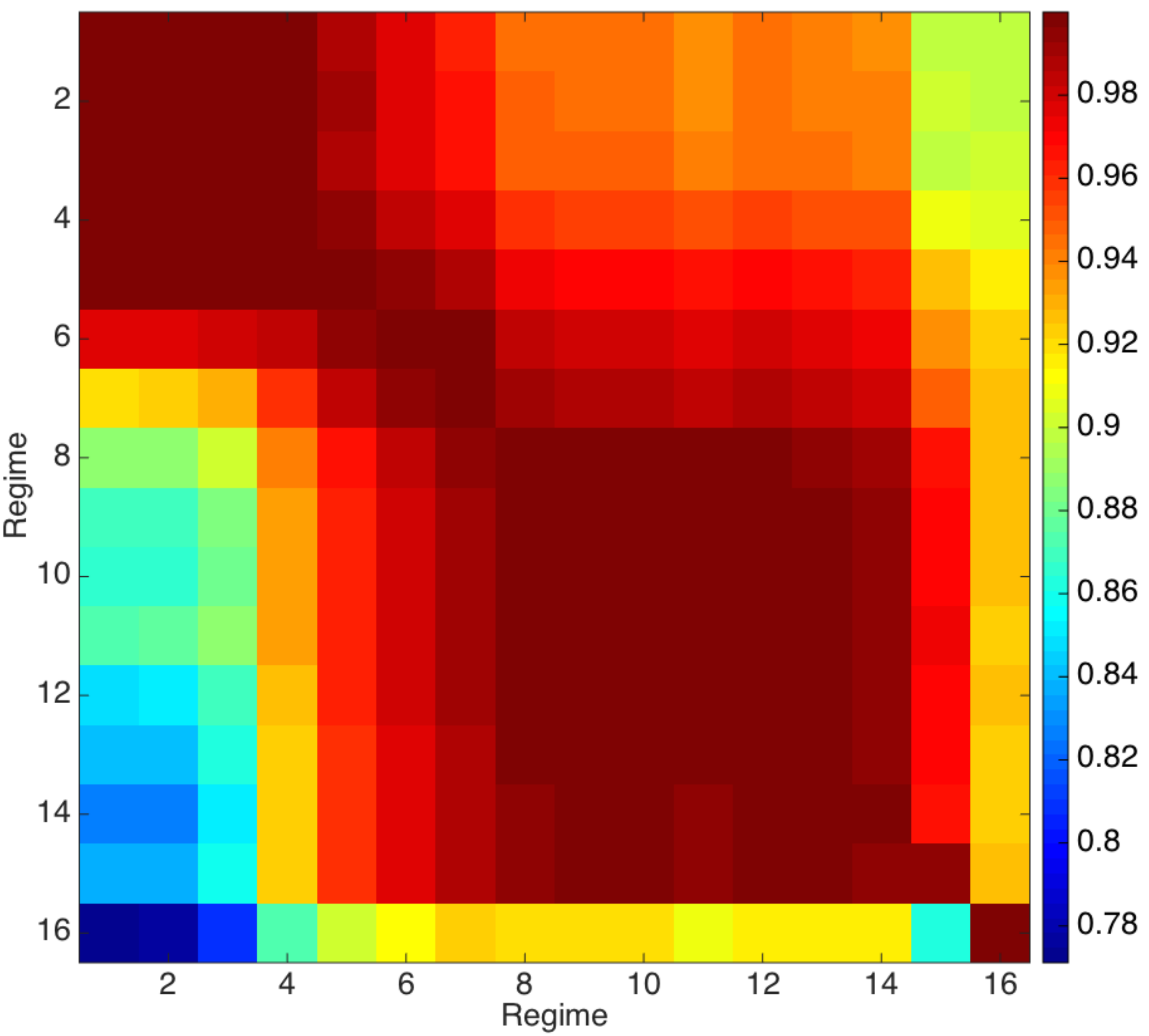}
	\subcaption{Projection onto full $n$-dimensional space, spanned by DMD basis functions, $P_i = \Phi \Phi_i^\dagger$.}	
	\end{subfigure}	
	\hspace{1cm}
	\begin{subfigure}[t]{0.45\textwidth}
	\includegraphics[width=8cm]{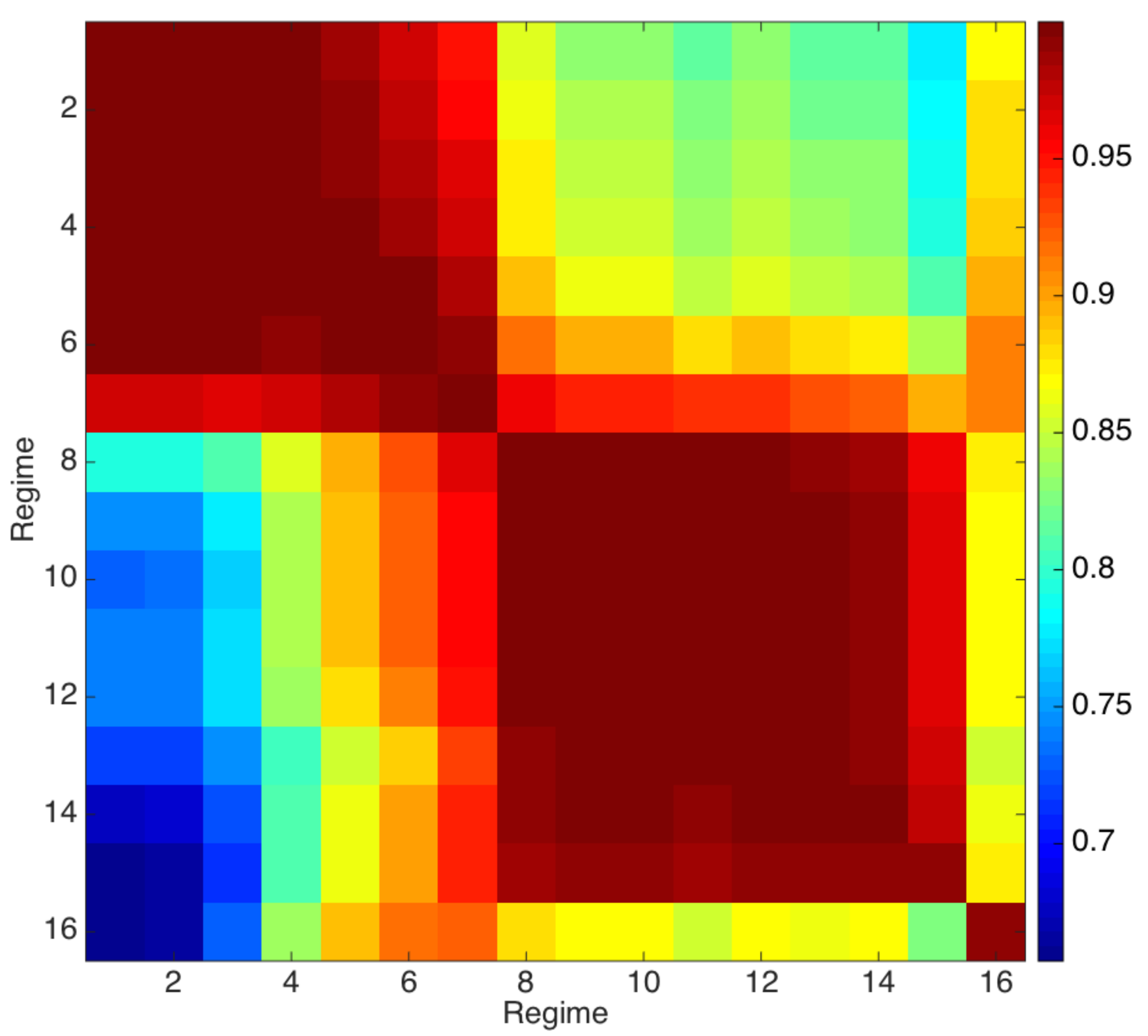}	
	\subcaption{Projection of boundary data $Y_i = CX_i$ with $P_i=\Theta_i\Theta_i^\dagger$ with $\Theta_i = C\Phi_i$.}	
	\end{subfigure}
	\caption{The data alignment measure $\kappa_{ij}$ from definition (\ref{def_kappa}).}
	\label{fig:DNS_dataalign}
\end{figure}

Figure \ref{fig:DNS_dataalign} shows this measure for the case of full
data and boundary data. The measure $\kappa_{ij}$, as defined in
equation (\ref{def_kappa}) indicates how much information is preserved
by projecting on the basis $\Phi_i$, through the projection $P_i$. As
shown in section \ref{section:DMD}, DMD modes are in general
non-orthonormal vectors that span the same subspace as singular
vectors given by the SVD step. Since we picked the top 20 singular
vectors to form the DMD subspaces $\Phi_i$, we see that the diagonal entries
for $\kappa_{ij}$ are mostly above 98\%. Additionally, neighboring
regimes share similar features; for instance, the projection of the
data from regime 3 onto the basis of regime 1 retains a high amount of
energy (measured in the Frobenius norm). As before, two groups of
regimes appear. For boundary data case, in analogy to the
considerations above for the measure $\gamma_{ij}$, the similarity
within the two clusters increases and the distinction among the two
clusters increases. Consequently, we expect some confusion within the
two clusters.

\subsection{Regime Classification}
\begin{figure}
	\centering
	
	\begin{subfigure}[t]{0.45\textwidth}
	\includegraphics[width = 8cm]{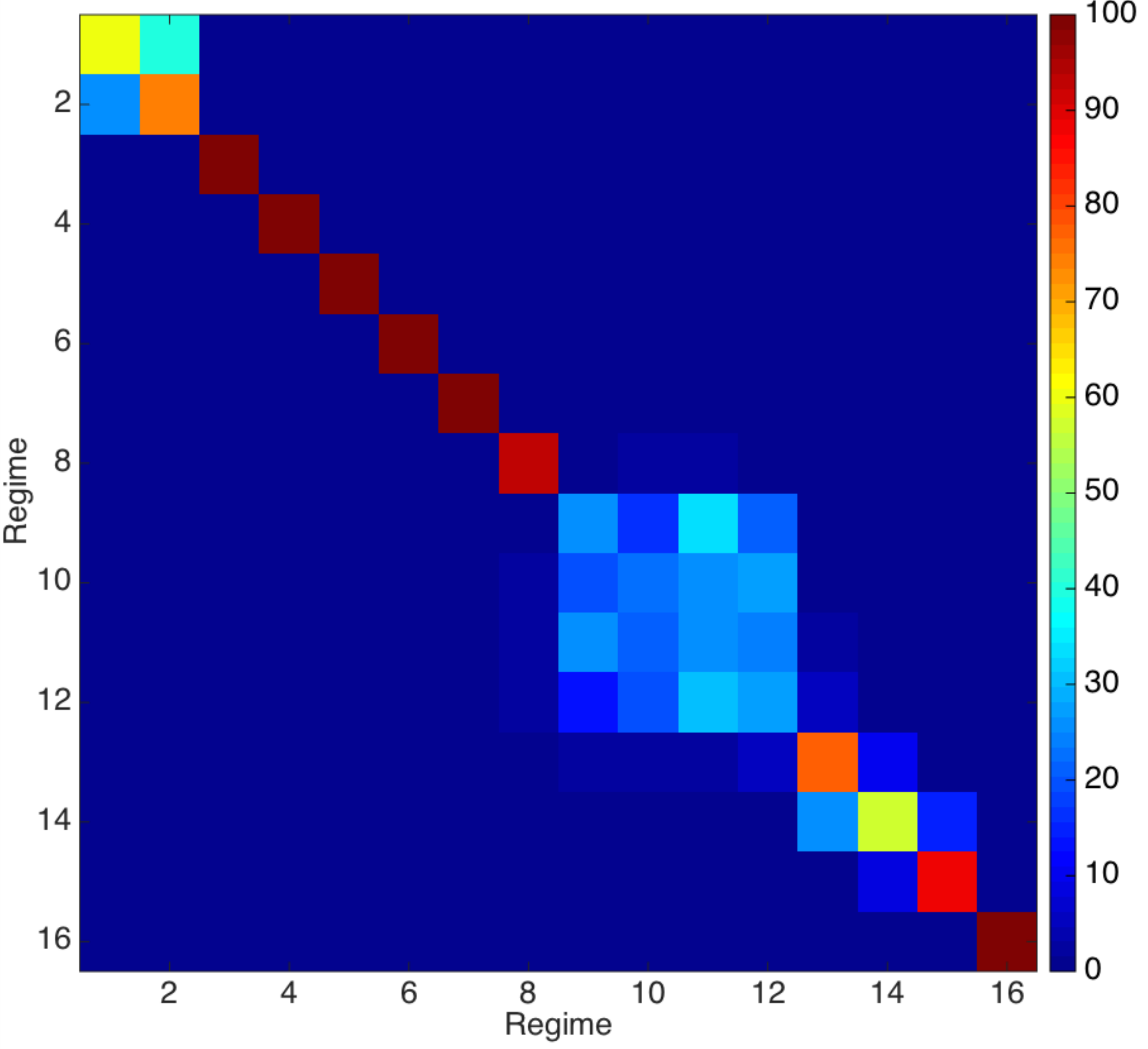}
	\subcaption{Flow regimes R1--R16.}	
	\end{subfigure}
\hspace{1cm}
	\begin{subfigure}[t]{0.45\textwidth}
	\includegraphics[width = 8cm]{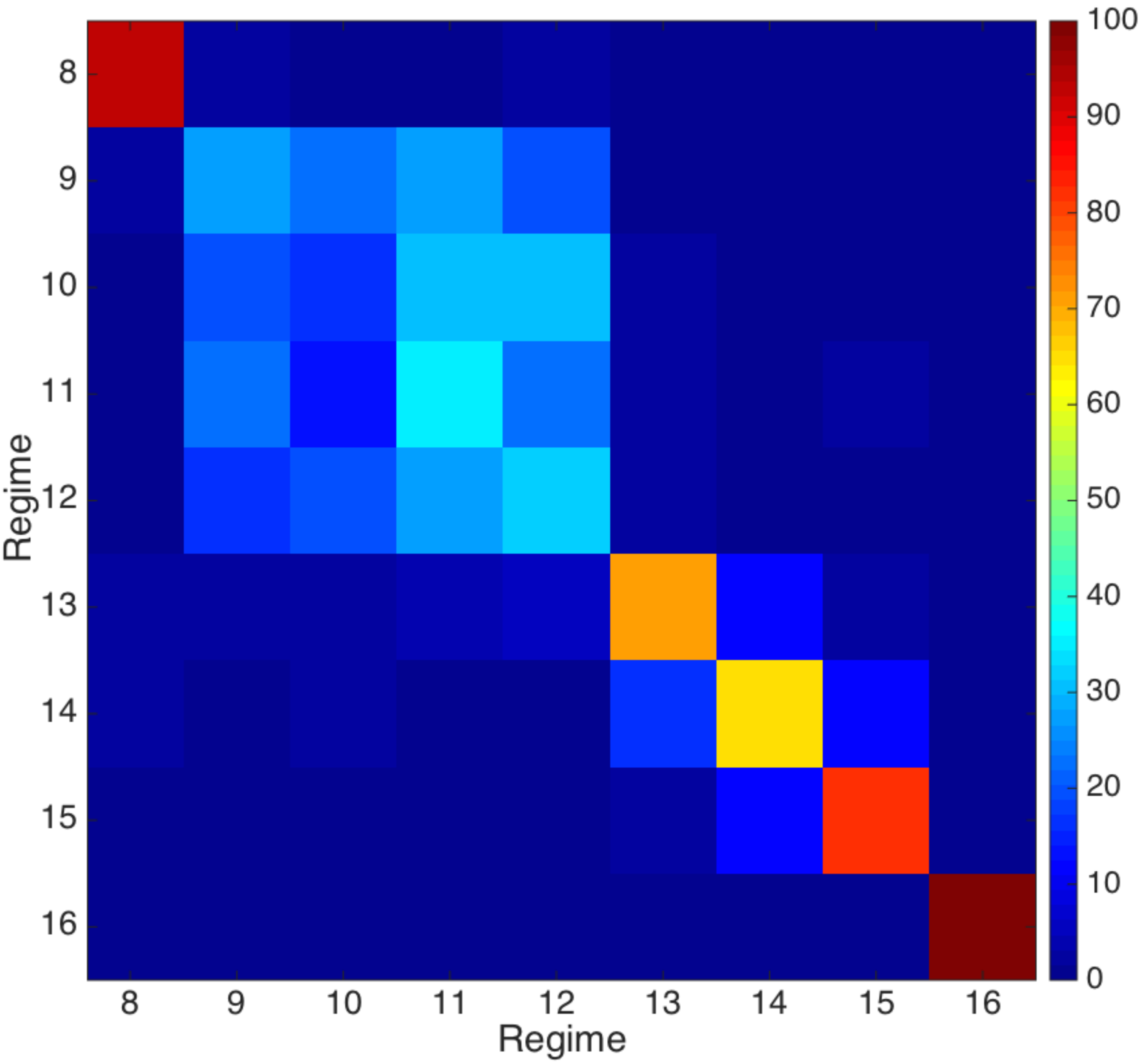}	
	\subcaption{Flow regimes R8--R16.}	
	\end{subfigure}
	\caption{Confusion matrices from 100 tests for every regime, indicating if the regime was picked correctly, or misclassified. In each test, we selected boundary sensors ($p_T =50$ for temperature and $p_v = 10$ for velocity), and used $j=3$ time measurements with augmented DMD.}
	\label{fig:DNS_conf}
\end{figure}

For the first example, we use $10$ velocity sensors, and $50$ temperature sensors. These sensors are placed on the boundary. The velocity is sensed one grid point away from the boundaries. The signal to noise ratio is set to $
20$dB, which corresponds to $10\%$ noise in $l_2$ sense. The DMD basis is augmented, i.e., we use the time evolution structure as described in Section \ref{secAugDMD} with $j=3$. For each regime, 100 tests are performed, where at each test, a snapshot from a given regime is picked and the best match to one of the 16 regimes is found by projection (\ref{argmax}). In Figure \ref{fig:DNS_conf}(a), the confusion matrix for the 16 regimes is plotted. The block with the highest confusion (worst identification) is between regimes R9--R12, corresponding to Rayleigh numbers $1.82\times10^8 - 2\times 10^8$, which is expected since the corresponding Rayleigh numbers are very close to each other. 
From the discussion in Section~\ref{section:bouss}, it is clear that the system has a steady state behavior for Rayleigh numbers corresponding to regimes R1 through R7. Thus, we further focus the analysis and discussion on Rayleigh numbers $Ra=10^8$ and higher, corresponding to R8--R16. The results, for another 100 tests for each regime, are presented in more detail in Figure \ref{fig:DNS_conf}(b).
\subsubsection{Robustness to out-of-sample data}
In the tests reported so far, all available data is used to generate the sparse library, and subsequently, the same data is used for classification. Here, training and testing datasets are separated, to investigate the robustness of the algorithm to unknown data. The testing data is taken from regimes R9, R10, and R11, which correspond to Rayleigh numbers between $1.82 \times10^8$ and $1.85\times10^8$. As noted in Section~\ref{section:bouss}, there is a bifurcation of the flow at $Ra \approx 1.82\times10^8$, which has been observed both experimentally, as well as numerically. We performed 600 independent tests, where at each test a flow snapshot (or a series of snapshots for the augmented sensing algorithm) is taken from the test data, and classified using six regimes R8, and R12--R16. One expects the classification of the testing data to match to regime R12, which is closest to the test regime. The signal-to-noise ratio is set to $20$dB, and $p_v=10$ flow sensors are used, together with $p_T=50$ temperature sensors, all placed on and near the boundary of the unit square. To be precise, the velocity sensors are placed slightly inside the domain, since the velocity is set zero at the boundaries. For better classification performance, and more robustness to noise, the DMD basis is augmented by two blocks, i.e., three time snapshots are taken for classification. 

\begin{table}[h!]
	\centering
	\begin{tabular}{c | c c c c c c}
		Regime 			&  R8 & R12 & R13 & R14 & R15 & R16  \\ \hline
				$Ra$ & $10^8$ & $2\times10^8$ & $4 \times10^8$ & $6\times10^8$ & $8 \times10^8$  & $10^9$\\ 
		Classification  & 6\%  & 87\% & 6\% & 0\% & 0\% & 0\% \\
	\end{tabular}
	\caption{Classification performance with out-of-sample sensor data (taken from regimes R9--R12) for six regimes at high Rayleigh numbers. The data is mainly classified as regime R12. }
	\label{tbl:DNS_unknwn_class} 
\end{table}

Table \ref{tbl:DNS_unknwn_class} reports the classification results. The sensing method is able to match the testing data to the (physically) correct flow patterns. In practice, this is important, since one does not expect the data to repeat itself in a given situation. Hence, the sensing mechanism needs to be able to match data to its ``closest" subspace in the library collection. \\
The regimes R9--R12 correspond to very close parameter values, and their solution subspaces are very similar to each other, as confirmed by looking at Figure \ref{fig:DNS_conf}, and the results of Table \ref{tbl:DNS_unknwn_class}. Hence, in the remaining analysis, we only consider the six regimes R8, and R12--R16, where R12 is the representative of regimes R9 to R12.

\subsubsection{Varying the number of sensors and using more data in time}
\begin{figure}[h]
	\centering
\includegraphics[width=5in]{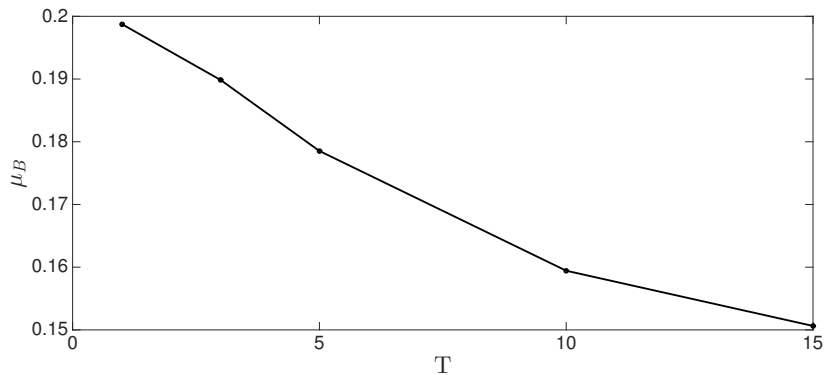}
 \caption{Block coherence measure $\mu_B$  as a function of time augmentation $T=j$ using only sensors on the boundary.}
 \label{fig:muB}
\end{figure}

\begin{figure}[h!]
\begin{center} 
\includegraphics[width=16cm]{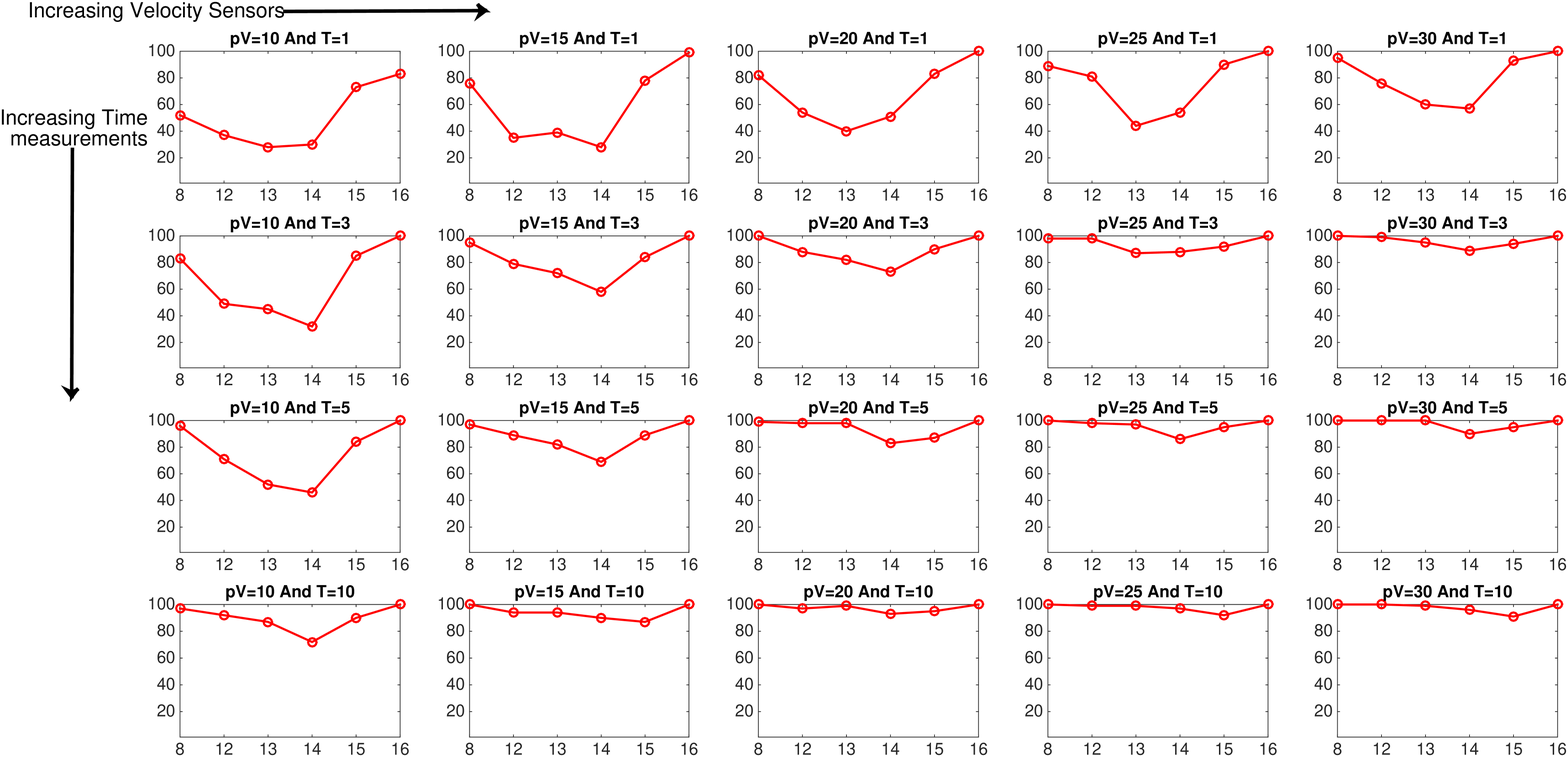}
\caption{Classification performance with $10\%$ sensor noise (SNR$=20$dB) for different number of near-boundary velocity sensors $p_v$, and time measurements $T=j$, with $p_T=10$ temperature sensors at the boundary. The performance increases significantly as number of time measurements are increased from 1 to 10, showing the efficacy of augmented DMD.}
\label{fig:result_matrix_SNR20}
\end{center}
\end{figure}

\begin{figure}[h!]
\begin{center} 
\includegraphics[width=16cm]{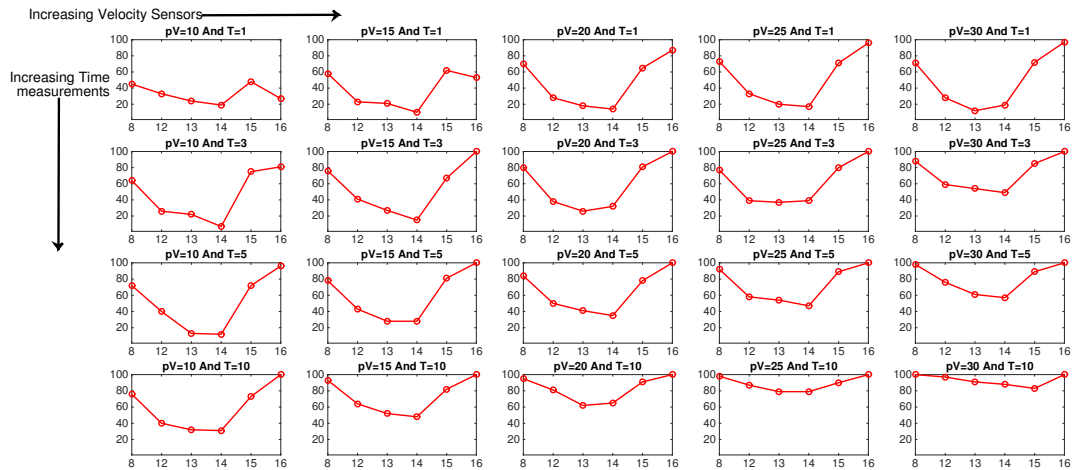}
\caption{Classification performance with $30\%$ sensor noise (SNR$=10$dB) for different number of near-boundary velocity sensors $p_v$, and time measurements $T=j$, with $p_T=10$ temperature sensors at the boundary. The performance increases significantly as number of time measurements are increased from 1 to 10, showing the efficacy of augmented DMD.}
\label{fig:result_matrix_SNR10}
\end{center}
\end{figure}

Next, we numerically study the improvement in classification performance that can be achieved by increasing the number of velocity sensors and time measurements, via the augmented DMD algorithm. In Figure \ref{fig:muB}, we plot the decay in block coherence measure $\mu_B$ with respect to the number of time measurements used in augmented DMD. In Figure \ref{fig:result_matrix_SNR20}, we plot the classification performance with $p_v=10,15,20,25,30$ velocity sensors, and 10 temperature sensors, all near the boundary. The number of measurements in time used in augmented DMD are $j=1,3,5,10$, and the sensor noise is $10\%$ ($SNR=20$dB). Figure \ref{fig:result_matrix_SNR10} shows the classification results with $30\%$ sensor noise ($SNR=10$ dB). It is evident from these figures that using the augmented DMD approach significantly improves the classification performance. While the value of the sub-coherence $\nu$ is still close to $0.99$, and hence the theoretical guarantee for accurate \emph{reconstruction} is highly conservative, the classification performance is much better. This improvement is reflected in the decay of $\mu_B$ as seen in Figure \ref{fig:muB}. While we have not focused on optimal sensor placement in this work, this connection between $\mu_B$ and classification performance can be used to formulate the problem of optimal sensor placement.

\clearpage

\subsection{Reconstruction in the Library}
The goal of this section is to show that the selection of the correct regime basis, i.e., the local model, is important for subsequent state estimation. To do so, we conduct the following simple numerical test for two different cases. We run the system in a specific regime $k^{\opt}$. Next, using $p_T=10$ temperature sensors and $p_v=30$ velocity sensors, we collect snapshots of the corresponding output measurement vector $y_{k^{\opt}}(t)\in\mathbb{R}^{70}$. This output vector is then projected onto bases of different regimes from the library $\mathcal{L}_{obs}$, and the resulting coefficients used to reconstruct the full $n$-dimensional state via the DMD library $\mathcal{L}$. By using the simple output-state map inversion 
$$
x_{i}^{\opt}(t)=\Phi_{i}\Theta_{i}^\dagger y_{k^\opt}(t), \quad i \in \{8,12,13,14,15,16\},
$$
we then obtain an estimate of the full state vector history. 

\begin{figure}[h!]
\begin{center} 
\includegraphics[width=.9\textwidth]{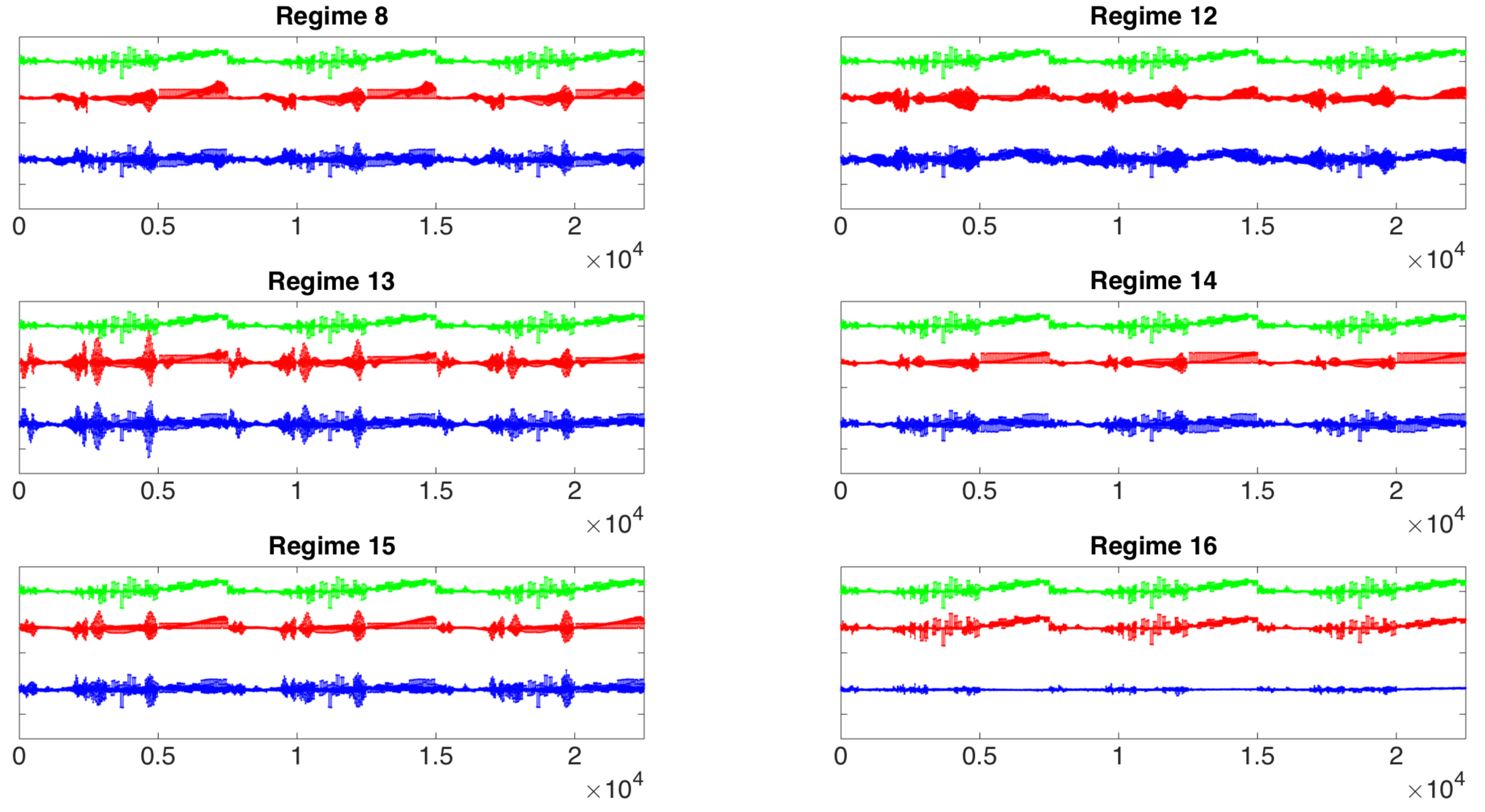}
\caption{\footnotesize{Reconstruction of three consecutive snapshots from a dynamic regime, corresponding to the parameter $Ra = 1\times 10^9$, in six different regimes, defined via $Ra = \{1\times 10^8, 2\times 10^8, 4\times 10^8, 6\times 10^8, 8\times 10^8, 1\times 10^9 \}.$ The x-axis is the index of the time augmented state vector. The original signal is shown in green, the reconstructed signal in red and the error is shown in blue. The corresponding relative errors in the Euclidean norm are $\{ .    0.5895, 0.8549, 0.7075, 0.5346, 0.5514, 0.1373\}$. The other parameters are:  $30\%$ Noise (SNR=$10$dB) , $p_v=30$ and $p_T=10$ }.} 
\label{fig:reconstr_DMD6}
\end{center}
\end{figure}

\begin{figure}[h!]
\begin{center} 
\includegraphics[width=.9\textwidth]{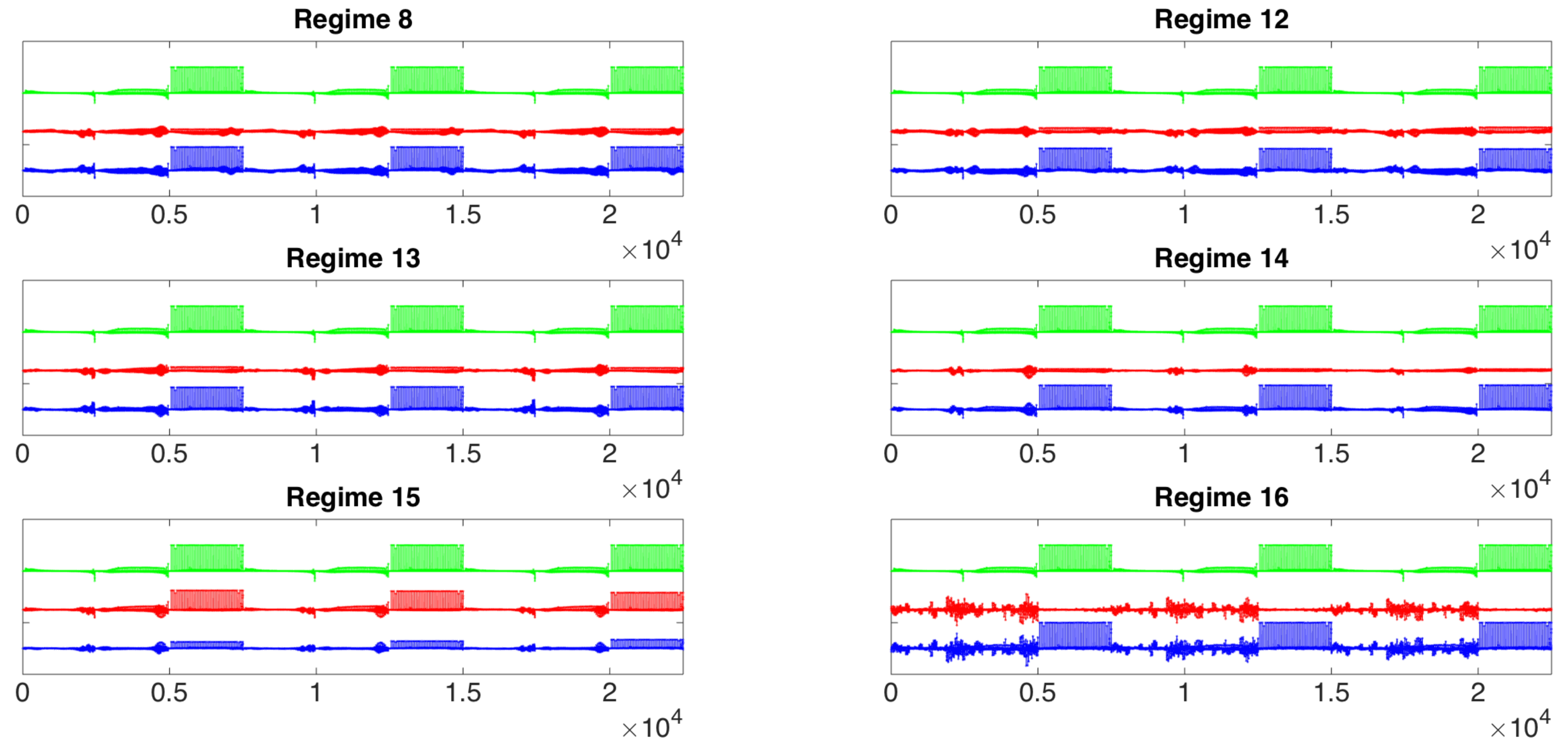}
\caption{\footnotesize{Reconstruction of three consecutive snapshots from a dynamic regime, corresponding to the parameter $Ra = 8\times 10^8$, in six different regimes, defined via $Ra = \{1\times 10^8, 2\times 10^8, 4\times 10^8, 6\times 10^8, 8\times 10^8, 1\times 10^9 \}.$ The x-axis is the index of the time augmented state vector. The original signal is shown in green, the reconstructed signal in red and the error is shown in blue. The corresponding relative errors in the Euclidean norm are $\{     1.0507,    0.9379 ,   0.9447  ,  0.9572  , 0.3806    ,1.1602\}$. The other parameters are:  $30\%$ Noise ($SNR=10$dB), $p_v=30$ and $p_T=10$ }
\label{fig:reconstr_DMD5}}
\end{center}
\end{figure}

In Figure \ref{fig:reconstr_DMD6}, we show the reconstruction results using three time measurements from regime R16 (i.e., $k^*=16$), and bases from the six regimes. In Figure \ref{fig:reconstr_DMD5}, we show the reconstruction results using three time measurements from regime R15 (i.e., $k^*=15$), and bases from the six regimes. Clearly, in both the test cases, the estimates using wrong regimes have significantly more error than estimates using the correct regime. 

Of course, this reconstruction is not a proper state reconstruction in the sense of state-space observation. Nevertheless, it shows the importance of selecting the right regime basis. After the correct regime basis has been selected, the corresponding low-dimensional model can be used in the framework of state observation. For instance, by using a Luenberger observer, one can properly estimate the dynamics of the full state vector over time. This part has not been studied in this paper, but will be the focus of our future work.


%% file: sectionconclusions.tex

\section{Conclusion}
We have introduced a framework for data-driven regime selection in parameter-dependent thermo-fluid systems. Most low-order modeling methods use
Galerkin projection-based dynamical reduced models of the infinite-dimensional fluid system. Hence the successful selection of dynamic regimes, and an
accurate reconstruction of the corresponding subspace is crucial to
the success of such models.

Our framework uses ideas from sparse sensing and exploits the
dimensionality reduction performed by reduced order models. In
particular, although our model can operate with a variety of dimensionality reduction methods, we use DMD to harness the ability of the Koopman
operator to capture the model dynamics. Thus, using a subspace
identification method, we can accurately identify the dynamic
regime with few sensors distributed near the boundary and
few time snapshots. 

Our framework uses off-line computation to construct a library of
regimes to be used for regime identification. This library comprises
of DMD eigenvectors and corresponding eigenvalues for each
regime, which capture the subspace in which the state lies under that
regime, as well as the dynamics of the state. The dynamical information is used to exploit multiple time measurements, which then increases robustness of the classification with respect to measurement noise and parametric uncertainty.

We numerically demonstrate  our approach using a DNS data set of a two
dimensional differentially heated cavity flow operating at different
Rayleigh numbers. The underlying PDE model is the Boussinesq equation,
which captures the dynamics of buoyancy driven flows well, when
temperature differences are small. The numerical results suggest that the proposed DMD-based classification method with augmented DMD basis is superior to using only a single-time measurement for classification.

Of course, our work is only a small step towards our goal, namely
low-order sensing and control models for complicated parametric
thermo-fluid systems. Parametric DMD~\cite{sayadi2015parametrized} is
one recent attempt to explicitly take into account the effect of
parameters in the DMD framework. State estimation has not been studied here, and will be the focus of our future report. Indeed, using the identified low-order model and
the correct regime it is possible to construct a stabilized low-order
model and correctly identify its parameters using data-driven learning
techniques~\cite{2015arXiv151001728B,benosmanbook}. Combining the DMD model with low-order Galerkin models is another promising
direction~\cite{noack14MORwithDMD}. These low dimension models can then be used for state observation, e.g., using Luenberger observers.

Regime construction is currently a manual process, understood
only for a few well-established systems. Automated optimal labeling
and sorting of data into different regimes, according to some measures
such as those discussed in this paper, would remove the need for
manual construction of regimes.

Another avenue of further research suggested by this work is the
connection between nonlinear systems and sparsity theory. While some
basic results in observability, controllability and state estimation
for linear systems with sparse states are already known
\cite{wakin10observability, sankaranarayanan2013compressive}, our
results suggest that this framework could be extended to nonlinear
systems via the DMD.
